\documentclass[preprint]{imsart}

\usepackage{
amsthm,amsmath,natbib,latexsym,graphicx,verbatim,amsfonts,subfigure}

\startlocaldefs

\newtheorem{Theorem}{Theorem}

\newtheorem{lemma}[Theorem]{Lemma}

\newtheorem{proposition}[Theorem]{Proposition}

\newcommand*{\cov}{\mathrm{cov~}}

\def\del{\delta}

\def\eps{\epsilon}
\def\kap{\kappa}

\def\sig{\sigma}

\def\bea{\begin{eqnarray*}}
\def\eea{\end{eqnarray*}}
\def\be{\begin{equation}}
\def\ee{\end{equation}}
\def\bean{\begin{eqnarray}}
\def\eean{\end{eqnarray}}
\def\barr{\begin{array}}
\def\earr{\end{array}}
\def\bdes{\begin{description}}
\def\edes{\end{description}}

\def\Bl{\Bigl}
\def\Br{\Bigr}

\def\til{\widetilde}
\def\hat{\widehat}

\def\Pr{\mathop{\rm I\!P}\nolimits}

\def\E{\mathbb{E}}
\def\N{\mathbb{N}}
\def\Q{\mathbb{Q}}
\def\R{\mathbb{R}}

\def\argmax{\mathop{\rm arg\,max}}
\def\argmin{\mathop{\rm arg\,min}}

\def\Bin{{\rm Bin}}

\def\supp{{\rm supp}}

\def\R{\mathbb{R}}
\def\E{\mathbb{E}}
\def\H{\mathbb{H}}
\def\N{\mathbb{N}}
\def\P{\mathbb{P}}
\def\Q{\mathbb{Q}}

\def\Bin{\mathrm{Bin}}

\def\AA{\mathcal{A}}
\def\BB{\mathcal{B}}
\def\CC{\mathcal{C}}
\def\DD{\mathcal{D}}
\def\FF{\mathcal{F}}
\def\GG{\mathcal{G}}
\def\HH{\mathcal{H}}
\def\JJ{\mathcal{J}}

\def\LL{\mathcal{L}}
\def\NN{\mathcal{N}}

\def\SS{\mathcal{S}}
\def\TT{\mathcal{T}}
\def\XX{\mathcal{X}}
\def\YY{\mathcal{Y}}

\def\del{\delta}
\def\eps{\epsilon}

\def\sig{\sigma}

\def\Pr{\mathbb{P}}

\def\argmin{\mathop{\rm argmin}}
\def\argmax{\mathop{\rm argmax}}

\def\bea{\begin{eqnarray*}}
\def\eea{\end{eqnarray*}}
\def\be{\begin{equation}}
\def\ee{\end{equation}}

\endlocaldefs

\begin{document}

\addtolength{\baselineskip}{0.2\baselineskip}
\begin{frontmatter}

\title{Optimal Calibration for Multiple Testing against Local Inhomogeneity in Higher Dimension
}
\runtitle{Optimal Caliration for Multiple Testing
}

\author{\fnms{Angelika} \snm{Rohde}\ead[label=e1]{angelika.rohde@math.uni-hamburg.de}}
\address{Universit\"at Hamburg\\
Bundesstra{\ss}e 55\\
20146 Hamburg\\
Germany\\
\printead{e1}}
\affiliation{Universit\"at Hamburg~~}
\runauthor{A. Rohde}

\begin{abstract} 
Based on two independent samples $X_1,...,X_m$ and $X_{m+1},...,X_n$ drawn from
multivariate distributions
with unknown Lebesgue densities $p$ and $q$ respectively, 
we propose an exact multiple test  in order to identify
simultaneously 
regions of significant deviations between $p$ and $q$.   
The construction is built from
randomized nearest-neighbor statistics. It  does not require any preliminary
information about the multivariate densities such as compact support, strict positivity or
smoothness and shape properties. The properly adjusted multiple testing procedure is shown to be sharp-optimal for typical arrangements of the observation values which appear with probability close to one. The proof relies
on a new coupling Bernstein type exponential inequality, reflecting the non-subgaussian tail behavior of a combinatorial process.  
For power investigation of the proposed method a reparametrized minimax set-up is introduced, reducing the composite hypothesis "$p=q$" to a simple one with the multivariate mixed density $(m/n)p+(1-m/n)q$ as infinite dimensional nuisance parameter. Within this framework, 
the test is shown to be spatially and sharply asymptotically
adaptive
with respect to uniform loss on isotropic H\"older classes. The exact minimax risk asymptotics are obtained in terms of solutions of the optimal recovery.
\end{abstract}



\end{frontmatter}

\section{Introduction}
\label{Introduction}
Given two independent multivariate iid samples 
$$
X_1,...,X_{m}\textrm{~~~and~~~}X_{m+1},...,X_{n}
$$
with corresponding Lebesgue densities $p$ and $q$ respectively, we are
interested in 
identifying simultaneously subregions of the densities support where
$p$ deviates significantly from $q$ at prespecified but arbitrarily chosen
level $\alpha\in (0,1)$. For this aim a multiple test of the composite hypothesis $H_0:\,
p=q$ versus $H_A:\, p\not=q$ is proposed, built from a suitable
combination of randomized nearest-neighbor
statistics. The procedure does not require any preliminary
information about the multivariate densities such as compact support, strict positivity or
smoothness and shape properties, and it is valid
for arbitrary finite sample sizes $m$ and $n-m$.
The hierarchical structure of
p-values for subsets of deviation between $p$ and $q$ provides insight into the local power of
nearest-neighbor classifiers, based
on the training set $\{X_1,...,X_n\}$. Thus our method is of interest in particular if
the classification error depends strongly on the value of the feature vector, related to recent literature on classification procedures by Belomestny and Spokoiny (2007).

There is an extensive amount on literature concerning two-sample problems. Most
of it is devoted to the one-dimensional case as there exists the simple but powerful
``quantile transformation'', allowing  for distribution-freeness under
the null hypothesis of several test statistics. Starting from the classical
univariate 
mean shift problem (see e.g{.} H\'ajek and $\check{\text{S}}$id\'ak 1967), 
 more flexible alternatives as stochastically larger or omnibus
 alternatives have been
 investigated for instance by  Behnen, Neuhaus and Ruymgaart (1983), Neuhaus (1982, 1987), Fan
 (1996), Janic-Wr\'oblewska and
Ledwina (2000), and Ducharme and Ledwina (2003). Our approach is different in that it aims at spatially adaptive and simultaneous identification of local rather than global deviations. In the above cited literature asymptotic power is discussed against single directional alternatives tending to zero at a prespecified rate, typically
 formulated by means of the densities $\tilde{p}$ and $\tilde{q}$
 corresponding to the transformed observations $\tilde{X}_i=H(X_i)$, where $H$
 denotes the
 mixed distribution function with density $h=(m/n)p+(1-m/n)q$. Note that the mapping $H$ coincides with the inverse quantile transformation under the null.

For power investigation of our procedure a specific two-sample minimax set-up is introduced. It is based on a reparametrization of $(p,q)$ to a couple $(\phi,h)$, reducing the composite hypothesis "$p=q$" to the simple one "$\phi\equiv 0$" with the multivariate mixed density $h$ as infinite dimensional nuisance parameter. The reparametrization conceptionally differs from the above described transformation for the univariate situation as it cannot rely on the inverse mixed distribution function. Nevertheless it leads under moderate additional assumptions in that case to the same notion of efficiency. 
In order to explore the power of our method, the alternative is assumed to be of the form
\begin{equation}
\Big\{ (p,q):\ (m/n)p+(1-m/n)q=h,\, \phi\in \FF,\, \Arrowvert \phi\Arrowvert\geq \delta\Big\}\label{Menge}
\end{equation}
for fixed but unknown $h$, some suitably chosen (semi-)norm $\Arrowvert .\Arrowvert$, a constant $\delta>0$  and a given smoothness class $\FF$.  For any 
$\alpha\in (0,1)$  the quality of a statistical level-$\alpha$-test $\psi$
is then quantified by its minimal power
$
\inf\E_{(p,q)}\psi,
$
where the infimum is running over all couples $(p,q)$ which are contained in the set (\ref{Menge}). It is a general problem that an optimal solution $\psi$ may depend on $\FF$ and $h$. Since the smoothness and shape of a potential difference $p-q$ are rarely known in practice, it is of interest to come up with a procedure which does not depend on these properties but is (almost) as good as if they were known, leading to the notion of minimax adaptive testing as introduced in Spokoiny (1996). Note that here we have however $h$ as an additional infinite dimensional nuisance parameter.

The problem of data-driven testing a simple hypothesis is further investigated for instance by Ingster (1987), Eubank and Hart (1992), Ledwina (1994), Ledwina and Kallenberg (1995), Fan (1996) and D\"umbgen and Spokoiny (2001) among others%
,  the two-sample context by Butucea and Tribouley (2006). The  idea in common is to combine a family of test statistics corresponding to different values of the smoothing parameters, respectively; see, for instance, Rufibach and Walther (2008) for a general criterion of multiscale inference. The closest in spirit to ours is the procedure developed in D\"umbgen and Spokoiny (2001) within the continuous time Gaussian white noise model and further explored by D\"umbgen (2002), D\"umbgen and Walther (2008) and Rohde (2008), all concerned with univariate problems. Walther (2010) treats the problem of spatial cluster analysis in two dimensions.

The paper is organized as follows. In the subsequent section, a multiple randomization test is intro\-duced,  built from a combination of suitably  standardized nearest-neighbor statistics. Its calibration relies on a new coupling exponential bound and an appropriate extension of the multiscale empirical process theory. 
Asymptotic power investigations and adaptivity properties are studied in Section 3, where the reparametrized minimax set-up is introduced. It is shown that our procedure is sharply asymptotically adaptive with respect to $\sup$-norm $\Arrowvert\cdot\Arrowvert$ on isotropic H\"older classes $\FF$, i.e. minimax efficient over a broad range of H\"older smoothness classes simultaneously. The application to local classification is discussed in Section 4. The one-dimensional situation is considered separately in Section 5 where an alternative approach based on local pooled order statistics is proposed. In that case the statistic does not depend on the observations explicitly but only on their order which in contrast to nearest-neighbor relations is invariant under the quantile transformation. 
Section 6.1 is concerned with a decoupling inequality and the coupling exponential bounds which are essential for our construction. Both results are of independent theoretical interest. All proofs and auxiliary results about empirical processes are deferred to Section 6.2 and Section~6.3.

\vspace{-1mm}
\section{Combining randomized nearest-neighbor statistics}\label{sec: nearest neighbors}
The procedure below is mainly designed for dimension $d\geq 2$. The
univariate case contains a few special features and is considered separately in Section \ref{sec: d=1}. 
Let $\underbar{X}:=(X_1,...,X_n)'$ and denote by $\XX_n$ the pooled set of observations. 
For any $0\leq k\leq n-1$, the $k$'th nearest-neighbor of
$X\in\XX_n$ with respect to the {\it Euclidean distance} is denoted by
$X^k$; we define $X^0:=X$. Note that the nearest-neighbors are unique a.s. The weighted labels are defined as follows

\vspace{-1mm}
$$
\Lambda(X)\ :=\begin{cases}
\frac{n}{m} & \ \text{if~$X$ is contained in the first sample} \\
-\frac{n}{n-m}& \ \text{otherwise}.
\end{cases}\vspace{-0.3mm}
$$ 
In order to judge about some possible deviation of $p$ from $q$ on a given set $B\in\BB^d$, a natural statistic to look at is a standardized version of $ \hat{\P}_n(B)-\hat{\Q}_n(B)$ or more sophisticated,

\vspace{-2.1mm}
$$
\int_B k_B(x)\Big(d\hat{\P}_n(x)-d\hat{\Q}_n(x)\Big)\vspace{-0.7mm}
$$
for some kernel $k_B$ supported by $B$, where $\hat{\P}_n$ and $\hat{\Q}_n$ denote the empirical measures corresponding to the first and second sample, respectively. Note that the statistic is not distribution-free, and in order to build up a multiple testing procedure several statistics corresponding to different sets $B$ have to be combined in a certain way.

\vspace{-2mm}
\subsection{Local nearest-neighbor statistics} 
Let $\psi:[0,\infty)\rightarrow\R$ denote any kernel of bounded total variation with $
\max_{x\in[0,\infty)}\arrowvert \psi(x)\arrowvert=\psi(0)=1
$ and $\psi(x)=0$ for $x>1$.
We introduce the local test statistics
\vspace{-0.8mm}
\begin{align*}
T_{jkn}\ :=&\ \frac{\sqrt{(m/n)(1-m/n)}}{\gamma_{jkn}}\frac{1}{\sqrt{n}}\sum_{i=0}^k\psi\bigg(\frac{\Arrowvert X_j-X_j^i\Arrowvert_2}{\Arrowvert X_j-X_j^k\Arrowvert_2}\bigg)\Lambda(X_j^i)\\
=&\ \frac{\sqrt{(m/n)(1-m/n)}}{\gamma_{jkn}}\sqrt{n}\int\psi\bigg(\frac{\Arrowvert X_j-x\Arrowvert_2}{\Arrowvert X_j-X_j^k\Arrowvert_2}\bigg)\Big(d\hat{\P}_n(x)-d\hat{\Q}_n(x)\Big),
\end{align*}
where 
$$
{\gamma_{jkn}}^2\ :=\ \frac{1}{n-1}\sum_{i=0}^{n-1}\Biggl[\psi\bigg(\frac{\Arrowvert X_j-X_j^i\Arrowvert_2}{\Arrowvert X_j-X_j^k\Arrowvert_2}\bigg)-\frac{1}{n}\sum_{l=0}^{n-1}\psi\bigg(\frac{\Arrowvert X_j-X_j^l\Arrowvert_2}{\Arrowvert X_j-X_j^k\Arrowvert_2}\bigg)\Biggr]^2.\\
$$
Every $T_{jkn}$ is some in a certain sense standardized weighted average of the nearest-neighbor's labels and its absolute value should tend to be large whenever $p$ is clearly larger or smaller than $q$ within the random Euclidean ball with center $X_j$ and radius $\Arrowvert X_j-X_j^k\Arrowvert_2$. 

\subsection{Adjustment for multiple testing} The idea is to build up a multiple  test, combining all possible local statistics $T_{jkn}$. The typical way is to consider the distribution of the supremum $\sup_{j,k}T_{jkn}$, see, e.g. Heckmann and Gijbels (2004). The problem is that the distribution is driven by small scales with a corresponding loss of power at larger scales, as there are many more small scales which contribute to the supremum. Here, we aim at a supremum type test statistic

\vspace{-4mm}
$$
T_n\ :=\ \sup_{0< k\leq n-1~}\sup_{1\leq j\leq n}\Big\{\arrowvert T_{jkn}\arrowvert -C_{jkn}\Big\},
\vspace{-1mm}
$$
where the constants $C_{jkn}$ are appropriately chosen correction terms (independent of the label vector $\Lambda$) for adjustment of multiple testing within every "scale" $k$ of $k$-nearest-neighbor statistics. These correction terms in the calibration aim to treat all the scales roughly equally. 
Although the distribution of $T_n$
under the null hypothesis depends on the unknown underlying distribution $p=q$,
the conditional distribution $\LL_0(T_n\arrowvert\XX_n)$ of the above statistic is invariant under permutation of
the components of the label vector $\underline{\Lambda}$. Here and subsequently, the index "$0$" indicates the null hypothesis, i.e. any couple $(p,q)$ with $p=q$. Precisely, 
let the random 
variable $\Pi$  be uniformly
distributed on the symmetric group $\SS_n$ of order $n$, independent of  $\underbar{X}$. Then
$
\LL_0\bigl({T}_{n}\big\arrowvert \XX_n\bigr)\ =\ \LL\bigl({T}_n\circ\Pi\big\arrowvert \XX_n\bigr),
$
where $\big(T_n\circ {\Pi}\big)(\underline{\Lambda}):=T_{n}\big(\Lambda_{\Pi_1},...,\Lambda_{\Pi_n}\big)$. Elementary calculation entails that
$$
\E\Bigl(T_{jkn}\circ{\Pi}\Big\arrowvert \XX_n
\Bigr)\ =\ 0\ \ \ \ \textrm{and}\ \ \ \ \ 
 \mathrm{Var}\Bigl(T_{jkn}\circ {\Pi}\Big\arrowvert
\XX_n\Bigr)\
=\ 1.
$$
Thus the null hypothesis is satisfied if, and only if, the hypothesis of
permutation invariance (or complete randomness)  conditional on $\XX_n$ is satisfied.

An adequate calibration of the randomized nearest-neighbor statistics, i.e. the choice of smallest possible constants $C_{jkn}$, requires both, an exact understanding of their tail behavior and their dependency structure. Note that the randomized nearest-neighbor statistics have a geometrically involved dependency structure. Even in case of the rectangular kernel $\psi$ it depends explicitly on the "random design" $\XX_n$ which incomplicates the sharp-optimal calibration for multiple testing compared to univariate problems, where the dependency of the single test statistics remains typically invariant under monotone transformation of the design points. Also, the optimal correction originally designed for Gaussian tails in D\"umbgen and Spokoiny (2001) does not carry over as only the subsequent Bernstein type exponential tail bound is available. 

\paragraph{\bf A coupling exponential inequality} 
Based on an explicit coupling, the following proposition extends and  tightens the exponential bounds
derived in Serfling\- (1974) for a combinatorial process in the present framework. If not stated otherwise, the random variable $\Pi$ is uniformly
distributed on $\SS_n$, independent of $\underbar{X}$.

\begin{proposition}\label{??}
Let $T_{jkn}$ be as introduced above and define 
$$
\delta(m,n)\ :=\ \bigg(\E\min\Bigl(\frac{S}{m},\frac{n-S}{n-m}\Bigr)\bigg)^{-1}
\ \ \ \textrm{with}\ \ \ 
S\sim\Bin\bigl(n, m/n\bigr).
$$ 
Then
$$
\P\Bigl(\big\arrowvert T_{jkn}\circ\Pi\big\arrowvert
   \, >\, \delta(m,n)\eta\Big \arrowvert\XX_n \Bigr) \leq\
  2\,\exp\left(-\,\frac{\eta^2/2}{1 + \eta \,n^{-1/2}\gamma_{jkn}^{-1}\,R_{\psi}(m,n)}\right),
$$
where
$$
R_{\psi}(m,n)\ :=\ \frac{2\Arrowvert
  \psi\Arrowvert_{\sup}}{3}\frac{\max(m,n-m)}{\sqrt{m(n-m)}}.
$$
\end{proposition} 

\paragraph{\sc Remark} The expression $\delta(m,n)$ is the 
payment for decoupling which  appears by replacing the tail probability of
 a hypergeometric ensemble by that of the
Binomial analogon. For details we refer to Section \ref{sec: decoupling}. In the typical case $0<\liminf_n (m/n)\leq \limsup_n (m/n)<1$ we obtain $\delta(m,n)=1+O(n^{-1/2})$. Compared to results obtained for weighted averages of standardized, independent Bernoullis, the above Bernstein type appears to be nearly optimal, i.e. subgaussian tail behavior (with leading constant $1/2$) is actually not present.

\medskip
Via inversion of the above exponential inequality, additive correction terms $C_{jkn}$ for adjustment of multiple testing are constructed. The next Theorem motivates our approach. The construction is designed for typical arrangements of the observation values which appear with probability close to one. To avoid technical expenditure, we restrict our attention to compactly supported densities. $d_w$ denotes the dual bounded Lipschitz metric (see, e.g. van der Vaart and Wellner 1996) which generates the topology of weak convergence. "$\rightarrow_{\P_n}$" refers to convergence in probability along the sequence of distributions $(\P_n)$. 

\begin{Theorem}\label{Thm: limit distribution} 
Define the test statistic
$$
T_n\ :=\ \sup_{\substack{1\leq j\leq n\\ 0< k\leq n-1}}\Big\{\,\arrowvert
  T_{jkn}\arrowvert\ -\ C_{jkn}\Big\}
$$
with 
$$
C_{jkn}\ :=\ 3\,R_{n}\gamma_{jkn}^{-1}\delta(m,n)\Gamma_{jkn}\ +\ \delta(m,n)\sqrt{2\,\Gamma_{jkn}},
$$
where $R_{n}=n^{-1/2}R_{\psi}(m,n)$ and
$\Gamma_{jkn}:=\log\big(1\big/{\gamma_{jkn}}^2\big)$. 
Assume that the sequence of mixed densities
$h_n:=(m/n)p_n+(1-m/n)q_n$ on $[0,1]^d$ is equicontinuous and uniformly bounded
away from zero, while $0<\liminf_n m/n\leq \limsup_n m/n<1$. Then the sequence  $\LL\bigl(T_n\circ \Pi\big\arrowvert \XX_n\bigr)$ of conditional distributions is tight in $\big(\P_n^{\otimes m}\otimes\Q_n^{\otimes(n-m)}\big)$-probability. Additionally, 
$$
d_w\Big( 
\LL\bigl(T_n\circ \Pi\big\arrowvert \XX_n\bigr),\, 
\LL(T_{\H_n})\Big)\ \longrightarrow_{\P_n^{\otimes m}\otimes\,\Q_n^{\otimes(n-m)}}\ 0,
$$
where
$$
  T_{\H_n}\ :=\ \sup_{\substack{t\in[0,1]^d,\\ 0<r\leq{\underset{x\in[0,1]^d}{\max}} \Arrowvert x-t\Arrowvert_2}
}\left\{\ \frac{\big\arrowvert\int_{[0,1]^d}
  {\phi}_{rt,n}(x)\,dW(x)\big\arrowvert}{\gamma_{rt,n}}\
  -\ \sqrt{2\log\bigl(1/{\gamma_{rt,n}}^2\bigr)}\ \right\}
$$
with $W$ a standard Brownian sheet in $[0,1]^d$, $
  {\gamma_{rt,n}} := \big(\int_{[0,1]^d} {\phi}_{rt,n}(x)^2d
  x\big)^{1/2}\ $ and
$$
\phi_{rt,n}(x)\ :=\ \bigg[\psi\Bigl(\frac{\Arrowvert x-t\Arrowvert_2}{r}\Bigr)\ - \ \int_{[0,1]^d}\psi\Bigl(\frac{\Arrowvert z-t\Arrowvert_2}{r}\Bigr)h_n(z)dz\bigg]\sqrt{h_n(x)}.
$$ 
\end{Theorem}

The extra-term $3\,R_{n}\gamma_{jkn}^{-1}\delta(m,n)\Gamma_{jkn}$ in the constant $C_{jkn}$ results from the exponential inequality in Proposition \ref{??} and can be viewed as an additional penalty for non-subgaussianity.  The theorem entails in particular that the sequence $\LL(T_n\circ\Pi\,\arrowvert \XX_n)$ is weakly approximated in probability by a tight sequence of {\it non-degenerate} distributions $\LL(T_{\H_n})$ which indicates that our corrections $C_{jkn}$ are appropriately defined and cannot be chosen essentially smaller.  
Note that the approximation $\LL(T_{\H_n})$ depends on the unknown mixed distribution even under the null hypothesis.

\subsection{The multiple rerandomization test} Let
$\kappa_{\alpha}(\underbar{X}):=\argmin_{C>0}\big\{ \P\big(T_n\circ\Pi\leq C\arrowvert\,
\XX_n\big)\geq 1-\alpha\big\}$ denote the generalized $(1-\alpha)$-quantile of
$\LL\big(T_n\circ\Pi\big\arrowvert \XX_n\big)$. Then we propose the
conditional test
$$
\phi_{\alpha}(\underbar{X})\ :=\ \begin{cases}
0& \textrm{if}\ T_n\leq\kappa_{\alpha}(\underbar{X})\\
1& \textrm{if}\ T_n >\kappa_{\alpha}(\underbar{X}).
\end{cases}
$$
Our method can be viewed as a multiple testing procedure. For a
given set of observations $\{X_1,...,X_n\}$, the corresponding test statistic
exceeds the $(1-\alpha)$-quantile if, and only if, the random set
$$
\DD_{\alpha}\ :=\ \Big\{ B_{X_j}\big(\Arrowvert
  X_j^k-X_j\Arrowvert_2\big)\Big\arrowvert\ 1\leq j\leq n,\ 0< k\leq n-1;\
  T_{jkn}(\underbar{X})\ >\ C_{jkn}(\underbar{X})\ +\
   \kappa_{\alpha}(\underbar{X})\Big\}
$$
is nonempty, where $B_t(r)$ denotes the Euclidean ball in $\R^d$ with center $t$ and radius $r$. Since the test is valid conditional on the set of observations,
we may conclude that $p$ deviates from $q$ at significance level $\alpha$ on
{\it every} Euclidean ball $B_t(r)\in\DD_{\alpha}$. In order to reduce the computational expenditure and to increase sensitivity on smaller scales, one may restrict one's attention to pairs $(j,k)$ for $k\leq m$ for some integer $m\in\{1,...,n-1\}$. Note the validity of the test does not require any assumption about the densities - even not Lebesgue continuity.

\medskip
Recently, Walther (2010) proposes a multiple test for cluster analysis in two dimensions based on a suitable combination of local $\log$-likelihood ratio statistics, evaluated on a fixed choice of axis-parallel rectangles. These statistics are not linear in $\hat{\mathbb{P}}_n$ and $\hat{\mathbb{Q}}_n$, respectively, but result in a subgaussian tail-behavior.

\section{Asymptotic power}\label{sec: efficiency} 

\subsection{Minimax-efficiency and spatial adaptivity -- local alternatives I}
In this section we show that the above introduced multiple testing procedure possesses  optimality properties in a certain minimax sense. Nonparametric comparison of different samples was recently investigated in the minimax approach by Butucea and Tribouley (2006), in a rate-adaptive way and of a different sense from our results here. We focus mainly on the considerably more involved problem of efficient adaptivity.
Let us first introduce some notation. For any set $J\subset [0,1]^d$ and function $f$ from $[0,1]^d\rightarrow \R$, $\Arrowvert f\Arrowvert_J:=\sup_{x\in J}\arrowvert f(x)\arrowvert$. For any convex $I\subset \R^d$ let $\HH_d(\beta,L;I)$ denote the isotropic H\"older smoothness class, which for $\beta\leq 1$ equals
$$
\HH_d(\beta,L;I)\ :=\ \Big\{ \phi: I\rightarrow\R:\, \big\arrowvert \phi(x)-\phi(y)\big\arrowvert\leq L\Arrowvert x-y\Arrowvert_2^{\beta}\Big\}. 
$$ 
Let $\lfloor \beta\rfloor$ denote the largest integer strictly smaller than $\beta$. For $\beta>1$, $\HH_d(\beta,L;I)$ consists of all functions $f:I\rightarrow \R$ that are $\lfloor\beta\rfloor$ times continuously differentiable such that the following property is satisfied: if $P_y^{(f)}$ denotes the Taylor polynomial of $f$ at the point $y\in I$ up to the $\lfloor\beta\rfloor$'th order, 
$$
\Big\arrowvert f(x)-P_y^{(f)}(x)\Big\arrowvert\ \leq\ L\Arrowvert x-y\Arrowvert_2^{\beta}\ \ \text{for all}\ x,y\in I.
$$
In particular the definition entails that $f\in\HH_d(\beta,L;\R^d)$ implies $f\circ U\in\HH_d(\beta,L;\R^d)$ for every orthonormal transformation $U\in\R^{d\times d}$.  For any pair of densities $(p,q)$ on $[0,1]^d$, let $h(m,n,p,q)$ denote the corresponding mixed density $(m/n)p+(1-m/n)q$. Fix a continuous density $h>0$ and define
$
\FF_{h}^{(m,n)}(\beta,L)
$ 
to be the set of pairs of densities such that
$$ 
\phi(m,n,p,q)\ :=\ \frac{p-q}{\sqrt{h(m,n,p,q)}}\, \in\, \HH_d\big(\beta,L;[0,1]^d\big)\ \ \text{and}\ \ h(m,n,p,q)=h.
$$

\paragraph{\bf Reparametrizing the composite hypothesis} With the notation above,
 $$
p=h\cdot\Big(1+(1-m/n)\,\phi\big/\sqrt{h}\,\Big)\ \  \text{and}\ \  q=h\cdot\Big(1-(m/n)\,\phi\big/\sqrt{h}\,\Big).
$$
Consequently "$p=q$" is equivalent to "$\phi\equiv 0$", and if  $(m/n)p+(1-m/n)q=h$ is kept fixed, the composite hypothesis "$p=q$" reduces  to the simple hypothesis "$\phi\equiv 0$". In order to develop a meaningful notion of minimax-efficiency for the two-sample problem we treat subsequently the mixed density $h=h(m,n,p,q)$ as fixed but 
unknown infinite dimensional nuisance parameter for testing the hypothesis 
$$
H_0:\ \phi= 0\textrm{~~~versus~~~}H_A:\ \phi\not= 0.
$$
Note that in case that $h$ is uniformly bounded away from zero and $p$ is close to $q$, $\phi$ coincides approximately with the difference $2\big(\sqrt{p}-\sqrt{q}\big)$, see also the explanation subsequent to Theorem \ref{Thm: lower bound}.

\paragraph{\sc Remark}
It is worth being noticed that the optimal statistic for testing $H_0$ against any fixed alternative $\phi$ equals the likelihood ratio statistic
$$
\frac{d\P_{(m,n,p,q)}}{d\P_{(m,n,h,h)}}(\underbar{X})\ =\ \prod_{i=1}^m\Big(1+(1-m/n)\frac{\phi}{\sqrt{h}}(X_i)\Big)\prod_{j=m+1}^n\Big(1- (m/n)\frac{\phi}{\sqrt{h}}(X_j)\Big),
$$
whose distribution still depends on $h$ under the null. Here and subsequently, the subscript $(m,n,p,q)$ indicates the distribution with density $\prod_{i=1}^mp\prod_{i=m+1}^nq$. The rational behind the repara\-metrization is to eliminate the dependency on the nuisance parameter $h$ in the expectation under the null of the first and second order term of the $\log$-likelihood expansion, resulting in asymptotic independence of $h$ for its distribution under the hypothesis for any local sequence $(\phi_n)$.

\medskip
The subsequent theorem is about the lower bound of hypothesis testing within the above defined classes of densities. 

\begin{Theorem}[Minimax lower bound]\label{Thm: lower bound}  Let
$$
\rho_{m,n}\ := \Big(\frac{n\log n}{m(n-m)}\Big)^{\beta/(2\beta +d)} \ \ \text{and define}\ \ 
c(\beta,L)\ :=\ \biggl(\frac{2\,d\,L^{d/\beta}}{(2\beta
  +d)\Arrowvert
  \gamma_{\beta}\Arrowvert_2^2}\biggr)^{\beta/(2\beta+d)},
$$
where $\gamma_{\beta}$ defines the solution to the optimal recovery problem (\ref{eq: recovery})  below. 
Assume that the sequence of mixed densities $(h_{n})$ on $[0,1]^d$ is equicontinuous and uniformly bounded away from zero. Then for any fixed $\delta>0$ and every nondegenerate rectangle $J\subset [0,1]^d$,
$$
\limsup_{n\rightarrow\infty}\inf_{\substack{(p,q)\in\overset{~}{\FF_{h_n}^{(m,n)}}(\beta,L):\\ \Arrowvert \phi\Arrowvert_J\geq
  (1-\delta)c(\beta,L)\rho_{m,n}}}\E_{(m,n,p,q)}\,\psi_n\ \leq \ \alpha
$$
for arbitrary tests
$\psi_n$ at significance level $\leq\alpha$.
\end{Theorem}
Note that $\psi_n$ may depend on $(\beta,L)$ and even on the nuisance parameter $h_n$ as already does the Neyman-Pearson test for testing $H_0$ against any one-point alternative. 

\medskip
We now turn to the investigation of the test introduced in Section \ref{sec: nearest neighbors}.
To motivate the choice of an optimal kernel for our test statistics and its relation to the optimal recovery problem, let us restrict our consideration to the Gaussian white noise context, leading in case of univariate H\"older continuous densities on $[0,1]$ with $\beta>1/2$ to locally asymptotically equivalent experiments 
$$
dX_{1n}(t)=p_n(t)\,dt+\frac{\sqrt{h_n(t)}}{\sqrt{m}}\,dW_1(t)\ \ \ \text{and}\ \ \  dX_{2n}(t)=q_n(t)\,dt+\frac{\sqrt{h_n(t)}}{\sqrt{(n-m)}}\,dW_2(t)
$$
for two independent Brownian motions $W_1$ and $W_2$ on the unit interval (Nussbaum 1996, Theorem 2.7 with $f_0=h_n$ and Remark 2.8). A multiscale statistic built from standardized differences of kernel estimates
$$
\frac{\sqrt{(m/n)(1-m/n)}}{\Arrowvert \psi\sqrt{h_n}\Arrowvert_2}\int \psi(t)\Big(dX_{1n}(t)-dX_{2n}(t)\Big)
$$
(which is actually not admissible since $h_n$ is unknown in general) 
then yields a distribution under the null close to ours in Theorem \ref{Thm: limit distribution}, up to the fact that our local integrals in dimension one are taken with respect to a Brownian bridge, reformulated to a Wiener process integrand by change of the kernel.
Concerning the optimization of $\psi$, the quantity to be maximized within this Gaussian white noise context appears to be the expectation of the single test statistics under the least favorable alternatives as their variances do not depend on the mean. In case $h_n\equiv 1$ this expression equals
\begin{equation*}
\inf_{\substack{{\phi}\in\HH_1(\beta,L;[0,1]):\\
\Arrowvert {\phi}\Arrowvert_{J}\geq \delta}}\,\frac{\int{\phi}(t)\,\psi(t)\,dt}{\Arrowvert\psi\Arrowvert_2},
\end{equation*}\nopagebreak
leading to the dual representation of the optimal recovery problem (see Donoho 1994a).

\paragraph{\bf The optimal recovery problem in higher dimension} In the framework of isotropic H\"older balls, the optimal recovery problem leads to the solution $\gamma=\gamma_{\beta}$ of the optimization problem
\begin{equation}\label{eq: recovery}
\text{Minimize $\Arrowvert \gamma\Arrowvert_2$ over all $\gamma\in\HH_d\big(\beta,1;\R^d\big)$ with $\gamma(0)\geq 1$.}
\end{equation}
The closedness of $\HH_d(\beta,1;\R^d)\cap\big\{\gamma: \R^d\rightarrow{\R}\big\arrowvert\, \gamma(0)\geq 1\big\}$ in $L_2$ entails that the solution exists, its convexity implies furthermore uniqueness whence by isotropy of the functional class $\HH_d(\beta,1;\R^d)$ it must be radially symmetric. In case $\beta\leq 1$, one easily verifies that 
$
\gamma_{\beta}(x)=\psi_{\beta}\big(\Arrowvert x\Arrowvert_2\big)=\big(1- \Arrowvert x\Arrowvert_2 ^{\beta}\big)_+$. In its generality, the optimal recovery problem in higher dimension has not yet been investigated. Considering the partial derivatives of $\gamma_{\beta}$ along the coordinate axes entails that $\psi_{\beta}$ is necessarily contained in $\HH_1(\beta,L;\R)$. However, the transfered optimization problem
\begin{equation}\label{eq: rec 2}
\text{minimize $\int\psi(r)^2 \arrowvert r\arrowvert^{d-1}dr\ $ over all   $\psi$  with  $\psi\big(\Arrowvert .\Arrowvert_2\big)\in\HH_d(\beta,1;\R)$ and $\psi(0)\geq 1$}
\end{equation}
does not coincide with the univariate optimal recovery problem due to the additional weighting by $\arrowvert r\arrowvert^{d-1}$ which comes into play by polar coordinate transformation. Whether the solution of (\ref{eq: rec 2}) for $\beta>1$ is compactly supported or not is still open. For the case of univariate densities, it is known that the solution of the optimal recovery problem has compact support for any $\beta>0$ (Leonov 1997)
, but an explicit solution in case $\beta>1$ is known for $\beta=2$ only. Concerning details and advice on its construction, see Donoho (1994b) and Leonov (1999). For dimension $d>1$, see Klemel\"a and Tsybakov (2001). 

\medskip
\noindent
The next Theorem is about the asymptotic power of the multiple test developed in Section \ref{sec: nearest neighbors}. We restrict our attention to compact rectangles of $(0,1)^d$ to avoid boundary effects. This restriction may be relaxed by the use of suitable boundary kernels, extending those of Lepski and Tsybakov (2000) for the univariate regression case to higher dimension.

\begin{Theorem}[Adaptivity and minimax efficiency]\label{thm: Effizienz2} Let $\phi_{n,\alpha}^*$ denote the multiple rerandomization test at significance level $\alpha$, based on the kernel  $\psi_{\beta}I\{\cdot\,\geq 0\}$ rescaled to $[0,1]$. In case of unbounded support of $\psi_{\beta}$, we may use a truncated solution $\psi_{\beta,K}=\psi_{\beta}I\{0\leq \cdot\leq K\} $. Let $0<\liminf_n m/n\leq\limsup_n m/n<1$. Assume that $(h_n)$ is equicontinuous and uniformly bounded away from zero. Then for any fixed $\delta>0$, there exists a $K>0$ such that 
$$
\liminf_{n\rightarrow\infty}
\inf_{\substack{(p,q)\in\FF_{h_n}^{(m,n)}(\beta,L):\\ \Arrowvert \phi\Arrowvert_J\geq
  (1+\delta)c(\beta,L)\rho_{m,n}}}\P_{(m,n,p,q)}\big(\phi_{n,\alpha}^*=1\big)\
=\ 1
$$
for any nondegenerate compact rectangle $J\subset (0,1)^d$. 
\end{Theorem}

In particular, the test is sharp-optimal adaptive with respect to the second H\"older parameter $L$. While in view of the results in Ingster (1987) the optimal rate of testing may be expected, some technical effort had to be done to propose a calibration achieving even sharp minimax-optimality. 

\paragraph{\sc Remark}It is worth being noticed that the procedure achieves the upper bound uniformly over a large class of possible mixed densities.
The intrinsic reason is that conditioning on $\XX_n$ is actually equivalent to conditioning on $\hat{\H}_n$, which indeed is a sufficient and complete statistic for the nuisance functional $\H_n$.

\paragraph{{\sc Remark} \rm (Sharp adaptivity with respect to $\beta$ and $L$)} Our construction, including the procedure especially designed for the one-dimensional situation, involves one kernel, shifted and rescaled depending on location and volume of the nearest-neighbor cluster under consideration. Due to the dependency of the optimal recovery solution $\gamma_{\beta}$ on $\beta$, the corresponding test statistic $T_n=T_n(\beta)$ achieves sharp adaptivity with respect to the second H\"older parameter $L$ only. Taking in addition the supremum 
$
T_n^*\ :=\ \sup_{\beta\in [\beta_0,\beta_1]}T_n(\beta)
$
over all kernels $\gamma_{\beta}$ within a compact range $[\beta_0,\beta_1]\subset (0,\infty)$, sharp adaptivity with respect to both H\"older parameters may be attained, provided that the above supremum statistic still defines a tight sequence (in probability), i.e. the corresponding sequence of $1-\alpha$-quantiles $\kappa_{\alpha}^*(\underbar{X})$ was stochastically bounded. Then the convergence 
\begin{align*}
\P_{(m,n,p_n,q_n)}\Big(T_n^* > \kappa_{\alpha}^*(\underbar{X})\Big)\ \geq\ \P_{(m,n,p_n,q_n)}\Big(T_{\hat{j}_n\hat{k}_nn}(\beta)-C_{\hat{j}_n\hat{k}_nn}(\beta)\, > \, \kappa_{\alpha}^*(\underbar{X})\Big)\ \rightarrow\ 1\ \ \text{as}\ n\rightarrow\infty
\end{align*}
for any random couple $(\hat{j}_n,\hat{k}_n)$ and any choice of $\beta$ 
could be be extracted from the proof of Theorem \ref{thm: Effizienz2}. At least for $\beta\in[\beta_0,1]$ this tightness may be deduced from the fact that the unimodal and symmetric $\psi_{\beta}$ depends continuously on $\beta$ in the sup-norm -- in particular $\LL\big(\big(\int\phi_{rt}^{(\beta)}(x)dW(x)\big)_{(t,r)}\big)$ as defined in Theorem 2 with $\psi=\psi_{\beta}$ depends continuously on $\beta$ in the weak topology. A general investigation especially for $\beta>1$ is beyond the scope of this article. 

\medskip
The next theorem shows however that our procedure simply based on the rectangular kernel is rate-adaptive with respect to both H\"older parameters 
$(\beta,L)$. 
Due to the fact that it combines locally all nearest-neighbor scales at the same time, it even adapts to inhomogeneous smoothness of $p-q$, i.e. achieves {\it spatial adaptivity}. 

\begin{Theorem}[Spatial rate-optimality]\label{thm: spatial adaptivity}
Let $\phi_{n,\alpha}^*$ denote the multiple rerandomization test based on the rectangular kernel. Assume that $0<\liminf_n m/n\leq\limsup_n m/n<1$. Then for any fixed $k\in\N$ and parameters  $(\beta_1,...,\beta_k,L_1,...,L_k)$, $K>0$ and any collection of disjoint compact rectangles $J_i\subset [0,1]^d$, $i=1,...,k$, there exist constants $d_i=d(\beta_i,L_i,K)$ with
$$
\liminf_{n\rightarrow\infty} \inf_{\substack{(p,q):\\
(p-q)_{\arrowvert J_i}\in\HH_d\big(\beta_i,L_i;J_i\big)\\ \Arrowvert p-q\Arrowvert_{J_i}\geq\, d_i\,\rho_{m,n}(\beta_i),\\ h(m,n,p,q)_{\arrowvert J_i}\geq\, K}} \P_{(m,n,p,q)}\bigg(J_i\cap\DD_{\alpha}(\XX_n)\not= \emptyset\ \forall \, i=1,...,k\bigg)\
=\ 1.
$$
\end{Theorem}

\subsection{The stylized type of locally constant alternatives on small and large scales -- local alternatives II}
The results from the previous paragraph deal with small scales of different (arbitrary) order depending on the smoothness classes under consideration. In particular, the minimax lower bound is concerned with scales tending to zero as $m,n\rightarrow\infty$, and it is not yet clear that there is no substantial loss at rather large scales. The size of possible deviation $\Arrowvert \phi\Arrowvert_{\sup}$ and the scale (here $\sim (\Arrowvert\phi_n\Arrowvert_{\sup}/L)^{1/\beta}$) are linked in a specific way depending on the smoothness class under consideration, because the smoothness assumptions do not allow for arbitrarily fast decay to zero. The next theorem is different in spirit. We do not focus on smoothness classes but on stylized situations with $\phi$ being lower bounded by a "plateau" of absolute value $c_n/\sqrt{n\delta_n^d}$ within a ball $B_x(\delta_n)$. With $\lambda$ denoting the Lebesgue measure on $\BB([0,1]^d)$, define
\begin{align*}
\JJ^{(m,n)}_+&(c,x,\delta)\, 
:=\, \bigg\{p,q \ \text{$\lambda$-densities on}\ [0,1]^d:\ 
  \phi(m,n,p,q)(z) \,\geq\, \frac{c}{\sqrt{n\delta^d}}\ \forall\, z\in B_x(\delta),\, 0< c\leq\sqrt{n\delta^d} \bigg\}, 
\end{align*}

\vspace{-4mm}
\begin{align*}
\JJ^{(m,n)}_-(c,x,\delta)\, 
:=\, \bigg\{p,q \ \text{$\lambda$-densities on}\ [0,1]^d:\ 
  \phi(m,n,p,q)(z) \,\leq\, \frac{-c}{\sqrt{n\delta^d}}\ \forall\, z\in B_x(\delta),\, 0<c\leq\sqrt{n\delta^d}\bigg\}
\end{align*}
and
$$
\GG^{(m,n)}(c,x,\delta)\ :=\ \JJ^{(m,n)}_+(c,x,\delta)\cup\JJ^{(m,n)}_-(c,x,\delta).
$$  
\begin{Theorem}\label{thm: local alternatives}
Assume that $0<\liminf_n m/n\leq\limsup_n m/n<1$.

\medskip
(i) 
If $\psi_n$ is any sequence of tests at significance level $\alpha\in(0,1)$, then
$$
\inf_{(p,q)\in\GG^{(m,n)}(c_n,x,\delta_n)}\E_{(m,n,p,q)}\,\psi_n\ \rightarrow \ 1
$$
implies that $n\delta_n^d\rightarrow\infty$ and $c_n\rightarrow\infty$.

\medskip
(ii) If $\psi_{n,\alpha}^*$ decribes the multiple rerandomization test based on the rectangular kernel at significance level $\alpha \in(0,1)$, 
$$
\inf_{\substack{(p,q)\in\GG^{(m,n)}(c_n,x,\delta_n)\\ h(m,n,p,q)\geq K>0}}\E_{(m,n,p,q)}\,\psi_{n,\alpha}^*\ \rightarrow \ 1,
$$
provided that $n\delta_n^d\rightarrow\infty$ and $\sqrt{\log(1/\delta_n)}/c_n\rightarrow 0$.
\end{Theorem}
In particular our test is also consitent against local alternatives of the type $\kappa_n\phi/\sqrt{n}$ for $\kappa_n\rightarrow\infty$. Comparing (i) and (ii) demonstrates that the adaptive search for the location of deviations costs an additional logarithm of its inverse scale. One may read out of the proof that the restriction for the sequence $(c_n)$ in (ii) can be slightly refined. 

\section{Application to classification}\label{sec: classification}
Suppose we are given an iid sample $(X_i,Y_i)$, $i=1,...,n$, where the marginal distribution of $X_i$ is assumed to be Lebesgue-continuous with density $h$ on $\R^d$, and $Y_i$ takes values in $\{0,1\}$ with
$$
\P\Big( Y_i=1\Big\arrowvert X_i=x\Big)\ =\ \rho(x).
$$
Then $M:=\sum_{i=1}^nY_i\ \sim\ \Bin\big(n,\lambda\big)$ with $\lambda:=\int\rho(x)h(x)dx$. 
Assuming $\lambda\in (0,1)$ to be known, the question of local classification is to identify simultaneously subregions in $\R^d$ where $\rho$ deviates significantly from $\lambda$ which results in local testing the hypothesis
$$
H_0:\, \rho\, =\, \lambda\ \ \ \text{versus}\ \ \ H_A:\, \rho\, \not = \, \lambda.
$$
Imitating our procedure introduced in Section \ref{sec: nearest neighbors}, we may combine suitably standardized local weight\-ed averages of labels, but the standardization differs due to the fact that the sum of (strictly) positive labels is random and not fixed, in particular $Y_1,...,Y_n$ are stochastically independent. Consequently, we may then rely the procedure on the classical Bernstein exponential inequality for weighted averages of standardized Bernoullis. Of course, the optimal separation constant for testing "$\rho=\lambda$" within some Euclidean ball $B_t(r)$ and its complement depends on the amount of observations in $B_t(r)$, whence analogously to the consideration above for the two-sample problem we may use the reparametrization of $(\rho, h)$ to $(\phi, h)$ with
$$
\phi\, :=\, \frac{\rho -\lambda}{\lambda(1-\lambda)}\sqrt{h}.
$$
The power optimality results carry over to the classification context with similar arguments as used in the proof of Theorem \ref{thm: Effizienz2}. We omit its explicit formulation at this point. 

\section{Distribution-freeness via quantile transformation -- the case  d=1}\label{sec: d=1}
The one-dimen\-sional situation allows for an alternative and more elegant approach based on order relations. 
For let $X_{(1)},..., X_{(n)}$ denote the order statistic built from the pooled sample and define for any $0\leq j<k\leq n$ the local test statistics
\begin{align*}
{U}_{jkn} :&=
\frac{\sqrt{(m/n)(1-m/n)}}{{\eta}_{jkn}}\frac{1}{\sqrt{n}}\sum_{i=j+1}^k\hspace{-2mm}\psi\Big(\frac{i}{k-j}\Big)\Lambda(X_{(i)}),
\end{align*}
where
$$
{\eta_{jkn}}^2\ :=\ \frac{1}{n-1}\sum_{i=1}^n\biggl(\psi\Big(\frac{i-j}{k-j}\Big)-\frac{1}{n}\sum_{l=1}^n\psi\Big(\frac{l-j}{k-j}\Big)\biggr)^2.
$$
Compared to the procedure described in the previous section, we omit the explicit dependence of the weights on the observed values. Note that in contrast to nearest-neighbor relations, the order remains invariant under quantile transformation, i.e. $\text{rank}(H_n(X_i))=\text{rank}(X_i)$, resulting in distribution-freeness of the corresponding multiscale statistic under the null. Suppose the null hypothesis is satisfied for some Lebesgue continuous distribution on the real line. Then conditional on the order statistics as well as unconditional, the label vector is uniformly distributed on the set
$$
\bigg\{\Lambda\in\big\{n/m, -n/(n-m)\big\}^n:\, \sum_{i=1}^n\Lambda_i^{-1}=0\bigg\}.
$$ 
The described test statistics are local versions of classical Wilcoxon rank sum statistics. We omit any further investigation as the calibration for multiple testing can be done analogously to that proved in Theorem \ref{Thm: limit distribution} -- but keep in mind  that the approximating Gaussian multiscale statistic under the null hypothesis will be independent of the nuisance  functional $\H_n$ due to the quantile transformation. Note that the use of typical mathematical tools for power investigation of rank statistics like Hoeffding's decomposition is getting involved because the kernel $\psi_{\beta}$ for $\beta\leq 1$ is not differentiable.

\newpage
\section{Proofs and further probabilistic results}
\subsection{Decoupling inequality and coupling exponential bounds}\label{sec: decoupling}
This section contains the coupling exponential bounds, i.e. in this context for weighted averages from a hypergeometric ensemble. 
Using a different technique, namely an explicit coupling construction, the subsequent proposition extends results of Hoeffding (1963) on decoupling of expectations of convex functions in the arithmetic mean of a sample without replacement. Whereas in the latter case decoupling with constant $1$ is actually correct, a simple counterexample for an ensemble of two elements already shows that the result does not extend to arbitrary weighted averages, and some payment for decoupling appears to be necessary. 

\begin{proposition}[Decoupling inequality]\label{Prop: decoupling}
Let $Z_1,Z_2,..., Z_n$ be iid with  
$$
\P\bigl(Z_i=1\bigr)\ =\
\frac{m}{n}\textrm{~~~and~~~}\P\bigl(Z_i=0\bigr)\ =\ 1-\frac{m}{n},\ \ \ \ 0<m<n.
$$
Let $a\in\R^n$ with $\sum_{i=1}^na_i=0$ and $\Psi:\R\rightarrow \R$ be
convex. Then 
$$
\E\biggl(\Psi\Bigl(\sum_{i=1}^na_iZ_i\Bigr)\bigg\arrowvert
\sum_{i=1}^nZ_i=m\biggr)\ \leq\ \E\,\Psi\Bigl(\delta(m,n)\sum_{i=1}^na_iZ_i\Bigr),
$$
with 
$$
\delta(m,n)^{-1}\ :=\ \E\min\Bigl(\frac{S}{m},\frac{n-S}{n-m}\Bigr),\ \ \
S\sim\Bin\Bigl(n, \frac{m}{n}\Bigr). 
$$
In particular,
$\delta(m,n)^{-1}=1+O(n^{-1/2})$ for $m/n\rightarrow \lambda\in (0,1)$.
\end{proposition}

\paragraph{{\sc Proof}}  Let $X$ be uniformly distributed on the set 
$$
\Big\{x\in\{0,1\}^n:\ \sum_{i=1}^nx_i\ =\ m\Big\}
$$
and let $S\sim\Bin(n, m/n)$ such that $X$ and $S$ are independent. Define
$$
M\ :=\ \big\{ i:\ X_i\ =\ 1\big\}.
$$
Conditional on $X$ and $S$, the random vector $Z\in\{0,1\}^n$ is constructed as
follows:

\medskip
\noindent
If $S>m$, let $Z_i=1$ for all $i\in M$ and let $(Z_i)_{i\in M^c}$ be uniformly
distributed on the set 
$$
\Big\{ z\in\{0,1\}^{M^c}:\sum_{i\in M^c}z_i=S-m\Big\}.
$$
For $S\leq m$, let $Z_i=0$ for all $i\in M^c$ and let $(Z_i)_{i\in M}$ be
uniformly distributed on
$$
\Big\{ z\in\{0,1\}^M:\sum_{i\in M}z_i=S\Big\}.
$$
Note that $Z_1,...,Z_n$ are iid $\Bin (1,m/n)$. Then
\begin{align*}
\E\,\Psi\Bigl(\sum_{i=1}^na_iZ_i\Bigr)\ &=\
\E\E\biggl(\Psi\Big(\sum_{i=1}^na_iZ_i\Bigr)\bigg\arrowvert X,S\biggr)\\
&\geq\ \E\,\Psi\biggl(\E\Bigl(\sum_{i=1}^na_iZ_i\Big\arrowvert X,S\Bigr)\biggr)\
\ \ \ \ \textrm{(Jensen inequality)}\\
&=\ \E\,\Psi\biggl(I\{S\leq m\}\frac{S}{m}\sum_{i\in M}a_i\ +\
I\{S>m\}\Big(\sum_{i\in M}a_i+\frac{S-m}{n-m}\sum_{i\in M^c}a_i\Big)\biggr)\\
&=\ \E\,\Psi\biggl(I\{S\leq m\}\frac{S}{m}\sum_{i\in M}a_i\ + \
I\{S>m\}\frac{n-S}{n-m}\sum_{i\in M}a_i\biggr)\ \ \ \ \ 
\text{\big(since~\,$\sum_{i=1}^n a_i\, =\, 0\,$\big)}\\
&=\
 \E\,\Psi\biggl(\min\Big(\frac{S}{m},\frac{n-S}{n-m}\Big)\sum_{i=1}^na_iX_i\biggr)\\
&=\
\E\E\,\bigg[\Psi\biggl(\min\Big(\frac{S}{m},\frac{n-S}{n-m}\Big)\sum_{i=1}^na_iX_i\bigg)\bigg\arrowvert
X\biggr]\\
&\geq\ \E\,\Psi\biggl(\E\Big\{\min\Big(\frac{S}{m},\frac{n-S}{n-m}\Big)\Big\}\sum_{i=1}^n
a_iX_i\biggr)\ \ \ \ \ \textrm{(Jensen inequality)}.
\end{align*}
Furthermore,
\begin{align*}
\E\,\min\Big(\frac{S}{m},\frac{n-S}{n-m}\Big)\ &=\
 1-\E\Big(\frac{(S-m)_-}{m}+\frac{(S-m)_+}{n-m}\Big)\\
&\geq\ 1-\E\Big(\frac{\arrowvert S-m\arrowvert}{\min(m,n-m)}\Bigr)\\
&\geq\ 1-\frac{\lambda(m,n)}{\sqrt{n}}
\end{align*}
with $\lambda(m,n):=\sqrt{m(n-m)}\big/\min(m,n-m)$, which is uniformly bounded
    for $m/n\rightarrow\lambda\in (0,1)$. 
$~$\hfill$\square$ 

\medskip
Using the decoupling above, the next proposition presents the exponential bounds for the combinatorial process which are essential for our construction. It implies Proposition \ref{??} and improves in particular exponential tail bounds for the hypergeometric distribution of Serfling (1974) in the coefficient in front of $\eta^2$ for $m/n$ close to zero or one, moderate $\eta$  and large $n$. Note that this coefficient is crucial for the efficiency of the testing procedure. The results may also be compared with the decoupling based exponential tail bounds in de la Pe$\tilde{\text{n}}$a (1994, 1999).

\begin{proposition}[Coupling exponential inequalities]\label{Prop: Bernstein}
Let $Z_1,...,Z_n$ be iid with
$$
\P\bigl(Z_i=1\bigr)\ =\
\frac{m}{n}\textrm{~~~and~~~}\P\bigl(Z_i=0\bigr)\ =\ 1-\frac{m}{n},\ \ 0<m<n.
$$
Let $\psi_1,...,\psi_n$ real valued numbers with $\bar{\psi}$ its arithmetic mean and denote
$$
{\gamma_{m,n}}^2\ :=\ \mathrm{Var}\biggl(\sum_{i=1}^n\psi_iZ_i\bigg\arrowvert\sum_{i=1}^nZ_i=m\biggr)\ =\ \frac{m(n-m)}{n(n-1)}\sum_{i=1}^n\big(\psi_i-\bar{\psi}\big)^2.
$$
Then in case of $\gamma_{m,n}\not=0$,
\begin{align*}
\P\left(\bigg\arrowvert\frac{1}{\gamma_{m,n}}\sum_{i=1}^n\psi_i\Big(Z_i-\frac{m}{n}\Big)\bigg\arrowvert
    >\delta(m,n)\eta\Bigg\arrowvert \sum_{i=1}^nZ_i=m\right)\ &\leq\ 
  2\,\exp\left(-\frac{\eta^2/2}{1\ +\ \eta\,R(\psi,m,n)}\right)\\
&\leq\ 2\,\exp\left(-\frac{3\eta}{2c(m,n)}+\frac{9}{2c(m,n)^2}\right),
\end{align*}
where
$$
R(\psi,m,n)\ :=\ \frac{\max_i\arrowvert
  \psi_i-\bar{\psi}\arrowvert}{3\,\gamma_{m,n}}\max\Big(\frac{m}{n},\, 1-\frac{m}{n}\Big)\ \ \ \ \ 
\textrm{and}\ \ \ \ \ c(m,n)\ :=\ \frac{\max(m, n-m)}{\sqrt{m(n-m)}}.
$$
\end{proposition}

\paragraph{\sc Proof} With 
$$
M\ :=\ \frac{\max_i\arrowvert\psi_i-\bar{\psi}\arrowvert}{\gamma_{m,n}}\max\Big(\frac{m}{n},\,1-\frac{m}{n}\Big)
$$
we obtain for any $t>0$
\begin{align}
\P\Bigg(&\frac{1}{\gamma_{m,n}}\sum_{i=1}^n\psi_i\Big(Z_i-\frac{m}{n}\Big)
   \, >\,\delta(m,n)\eta\Bigg\arrowvert \sum_{i=1}^nZ_i=m\Bigg)\nonumber\\ 
&=\ \P\left(\frac{1}{\gamma_{m,n}}\sum_{i=1}^n\bigl(\psi_i-\bar{\psi}\bigr)\Big(Z_i-\frac{m}{n}\Big)\,
    >\,\delta(m,n)\eta\Bigg\arrowvert \sum_{i=1}^nZ_i=m\right)\nonumber\\
&\leq\ \exp\Big(-t\,\frac{\eta}{M}\Big)\,\E\biggl\{\exp\biggl(\frac{t\,\delta(m,n)^{-1}}{M\,\gamma_{m,n}}\sum_{i=1}^n\bigl(\psi_i-\bar{\psi}\bigr)\Bigl(Z_i-\frac{m}{n}\Bigr)\bigg)
     \bigg\arrowvert \sum_{i=1}^nZ_i=m\biggr\}\nonumber\\
&\leq\ \exp\Big(-t\,\frac{\eta}{M}\Big)\,
\E\exp\biggl(\frac{t}{M\,\gamma_{m,n}}\sum_{i=1}^n\bigl(\psi_i-\bar{\psi}\bigr)\Bigl(Z_i-\frac{m}{n}\Bigr)\biggr)\
\ \ \textrm{(Proposition \ref{Prop: decoupling})}\nonumber\\
&\leq\ \exp\biggl(\frac{1}{M^2}\big(e^t-1-t\big)\ -\ t\,\frac{\eta}{M}\biggr)\label{Bennett}, 
\end{align}
whereby the last inequality follows from the fact that for any random variable
$Y$ with $\arrowvert Y\arrowvert \leq 1$, $\E Y=0$ and
$\mathrm{Var}(Y)=\sigma^2$,
$$
\E\exp(tY)\ \leq\ 1+\sigma^2(e^t-1-t)\ \leq\ \exp\Big(\sigma^2(e^t-1-t)\Big).
$$
Elementary algebra shows that (\ref{Bennett}) is
minimized with the choice 
$
t := \log\big(1+ \eta M\big),
$
which yields first a Bennett (1962) exponential bound and because of $(1+x)\log(1+x)-x\geq (x^2/2)/(1+x/3)$ consequently the Bernstein type
$$
\P\Bigg(\frac{1}{\gamma_{m,n}}\sum_{i=1}^n\psi_i\Big(Z_i-\frac{m}{n}\Big)
   \, >\,\delta(m,n)\eta\Bigg\arrowvert \sum_{i=1}^nZ_i=m\Bigg)\ \leq\ \exp\biggl(-\frac{\eta^2/2}{1+\eta\,M/3}\biggr).
$$  
A symmetry argument provides the same bound for $\psi_i$ replaced by
$-\psi_i$, which completes the proof of the first inequality. Using that $\gamma_{m,n}\geq \sqrt{(m/n)(1-m/n)}\max_i\arrowvert\psi_i-\bar{\psi}\arrowvert$,  
we obtain 
the second
asserted inequality from
\begin{align*}
\frac{\eta^2/2}{1+\eta\,M/3}\  &\geq\ \frac{\eta^2/2}{1+\eta\,c(m,n)/3}\\
 &= \ \frac{\eta}{2c(m,n)/3}\ -\ \frac{\eta}{2c(m,n)/3(1+\eta\, c(m,n)/3)}\\
 &\geq \ \frac{\eta}{2c(m,n)/3}\ -\ \frac{1}{2c(m,n)^2/9}.
\end{align*}

\vspace{-8mm}
$~$\hfill$\square$

\subsection{Auxiliary results about empirical processes}
This section collects results in the context of empirical processes which are essential for the next section. For any totally-bounded pseudo-metric space $(\TT,\rho)$, we define the covering number
$$
N\big(\varepsilon,\TT,\rho\big)\ :=\ \min\Big\{\sharp \TT_0:\, \TT_0\subset\TT,\inf_{t_0\in\TT_0}\rho(t,t_0)\leq\varepsilon\text{~for all~}t\in\TT\Big\}.
$$
Let $\BB(\TT)$ denote the Borel-$\sigma$-field on $\TT$ induced by the pseudo-metric $\rho$ (which induces a topology in the usual sense, although without the Hausdorff-property if it is not a metric) and let $\mathcal{F}\subset [0,1]^{\mathcal{T}}$ be a family of measurable functions. For any probability measure $P$ on $\BB(\TT)$, consider the pseudo-distance $d_P(f,g) := \int \arrowvert f-g\arrowvert \, dP$ for $f, g \in \mathcal{F}$. Then for any $u>0$, the uniform covering numbers of $\mathcal{F}$ are defined as
$
	\mathcal{N}(u,\mathcal{F}):=\sup_{P}N(u,\mathcal{F}, d_P)
$,
 where the supremum is running over all probability measures $P$ on $\BB(\mathcal{T})$.

\begin{Theorem}{\rm (D\"umbgen and Walther (2008, technical report))}
\label{Chaining}
Let $Z = (Z(t))_{t \in \TT}$ be a stochastic process on a totally bounded pseudo-metric space $(\TT, \rho)$. Let $K$ be some positive constant, and for $\del > 0$ let $G(\cdot, \del)$ a nondecreasing function on $[0, \infty)$ such that for all $\eta \ge 0$ and $s, t \in \TT$,
\be
	\label{Chaining.1}
	\Pr \Bl\{ \frac{|Z(s) - Z(t)|}{\rho(s,t)} > G(\eta, \del) \Br\} 
	\ \le \ K \exp(-\eta) 
	\quad \mbox{if } \rho(s,t) \ge \del.
\ee
Then for arbitrary $\del > 0$ and $a \ge 1$,
$$
	\Pr \Bigl\{ |Z(s) - Z(t)| \ge 12 J(\rho(s,t),a)
		\mbox{ for some } s,t\in\TT_* \mbox{ with } \rho(s,t) \le \del \Bigr\} 
	\ \le \ \frac{K \del}{2a} \, ,
$$
where $\TT_*$ is a dense subset of $\TT$, and
\bea
	J(\eps,a) 
	& := & \int_0^\eps G(\log(a D(u)^2 / u), u) \, du, \\
	D(u) = D(u,\TT,\rho)
	& := & \max \Bl\{ \#\TT_o : \TT_o \subset \TT, \rho(s, t) > u \mbox{ for different } 
		s, t \in \TT_o \Br\}.
\eea
\end{Theorem}

\noindent
\textbf{Remark.} Suppose that $G(\eta, \del) = \til{q} \, \eta^q$ for some constants $\til q, q > 0$. In addition let
$
	D(u) \ \le \ A u^{-B}	\quad\mbox{for } 0 < u \le 1
$
with constants $A \ge 1$ and $B > 0$. Then elementary calculations show that for $0 < \eps \le 1$ and $a \ge 1$,
$
	J(\eps,a) \ \le \ C \, \eps \, \log(e/\eps)^q
$
with $C = \til{q} \, \max \bigl( 1+2B, \log(a A^2) \bigr)^q \int_0^1 \log(e/z)^q \, dz$.

\medskip
For the proof of Theorem \ref{Thm: limit distribution} the subsequent
extension of the Chaining Lemma VII.9 in Pollard (1984) and Theorem 8 in the technical report to D\"umbgen and Walther (2008) will be used. It complements in particular the existing multiscale theory by a uniform tightness result and to a situation where only a sufficiently sharp {\it uniform stochastic} bound on local covering numbers is available, which typically involves additional logarithmic terms. The situation arises for example in the multivariate random design case where a non-stochastic bound obtained via uniform covering numbers and VC-theory may be too rough. 

\begin{Theorem}[Chaining]
\label{Levy}
Let $(Y_n)_{n\in\N}$ be a sequence of random variables such that $Y_n$ takes values in some polish space $\YY_n$. For any $y_n\in\YY_n$, let $(Z_n(t;y_n))_{t\in\TT_{y_n}}$ be a stochastic process on some countable, metric space $\big(\TT_{y_n},\rho_n(.,.;y_n)\big)$, where $\rho_n(.,.;y_n)\leq 1$. 
Suppose that the following conditions are satisfied:

\medskip
\noindent
\textbf{(i)} \ There are measurable functions $\sig_n(.;Y_n) : \TT_{Y_n} \to (0,1]$ and $G_n(.,\delta) : [0,\infty) \to [0,\infty)$ such that for arbitrary $s,t \in \TT_{Y_n}$, $\eta \ge 0$ and $\del > 0$,
$$
	\P \Bigl( |Z_n(t,Y_n)| \geq \sig_n(t;Y_n)  G_n(\eta,\del)\Big\arrowvert Y_n \Bigr) 
	\ \le \ 2 \exp(- \eta)	\quad\mbox{if } \sig_n(t;Y_n) \geq \del ,
$$
$$
	\sup_{s,t\in\TT_{Y_n}}\frac{\arrowvert
          \sigma_n(t;Y_n)-\sigma_n(s;Y_n)\arrowvert}{\rho_n(s,t;Y_n)}\ \leq\
        C<\infty \       
	\mbox{for some constant $C>0$,}
$$
$$	\big\{t\in\TT_{Y_n}:\sig_n(t;Y_n)\ge\del\big\} \mbox{ is compact,\ \  and } 
	G_o := \sup_{n\in\N}\sup_{\eta \ge 0, 0 < \del \le 1} \frac{G_n(\eta,\del)}{1 + \eta}
	\ < \ \infty \, .
$$
\textbf{(ii)} \ There exists a sequence  $(\CC_n)_{n\in\N}$ of measurable sets and positive constants $A,B,W,\alpha$ such that
$$
	N \Bigl( u\del, \{t\in\TT_{Y_n} : \sig_n(t;Y_n) \leq \del\}, \rho_n(.,.;Y_n) \Bigr) 
	\ \leq \ A u^{-B} \del^{-W}\log\big(e/(u\delta)\big)^{\alpha}	\quad\mbox{for }
        u,\ \del \in (0,1] 
$$
whenever $Y_n\in\CC_n$.

\medskip
\noindent
For constants $q,Q > 0$ define 
$$
	\AA_n(\del,q,Q;Y_n) \ := \ \bigg\{ \sup_{s,t \in \TT_{Y_n} \, : \, \rho_n(s,t;Y_n) \le \del} \, 
		\frac{|Z_n(s;Y_n) - Z_n(t;Y_n)|}{\rho_n(s,t;Y_n) \log(e/\rho_n(s,t;Y_n))^q}\ \le\ Q \bigg\}.
$$

\noindent Then there exists a constant $C = C(G_o,A,B,W,\alpha, q,Q) > 0$ such that for $0 < \del \le 1$ 
\begin{align*}
	\P\Biggl(
        \frac{|Z_n(t;Y_n)|}{\sig_n(t;Y_n)} \, \le \,  G_n \biggl( W
        \log\big(1/\sig_n(t;Y_n)\big)\ 
+& \ C \log\log\big(e/\sig_n(t;Y_n)\big),\, \sig_n(t;Y_n) \biggr) \\
		+\ C  \log(e/\sig_n&(t;Y_n))^{-1} \mbox{ on } \big\{t :
                \sig_n(t;Y_n) \le \del\big\}\Bigg\arrowvert Y_n\Biggr)
\end{align*}
is at least 
$
\P\Big(\AA_n(2\del,q,Q;Y_n)\Big\arrowvert Y_n\Big) - C \log(e/\del)^{-1}
$
whenever $Y_n\in\CC_n$. 

If in particular $\P^{Y_n}(\CC_n)\rightarrow 1$ and $\lim_{\delta\searrow 0}\inf_n\P\Big(\AA_n(\del,q,Q;Y_n)\Big\arrowvert Y_n\Big)=1$ a.s., then the sequence
\begin{align*}
\LL\Bigg(\sup_{t\in\TT_{Y_n}}\Bigg\{\frac{|Z_n(t;Y_n)|}{\sig_n(t;Y_n)}& \, - \,  G_n \biggl( W
        \log\big(1/\sig_n(t;Y_n)\big)
 \ + \ C \log\log\big(e/\sig_n(t;Y_n)\big),\, \sig_n(t;Y_n) \biggr) \Bigg\}
	 \Bigg\arrowvert Y_n\Biggr)
\end{align*}
is tight in $\big(\P^{Y_n}\big)$-probability, provided that $\inf_n\sup_{t\in\TT_{Y_n}}\sigma_n(t;Y_n)>0$ a.s.

\end{Theorem}

\paragraph{\sc Remark} Note that in case of $G(\eta,\del) = (\kap \eta)^{1/\kap}$ with $\kap > 1$,
\bea
	\lefteqn{ G \Bl( W \log(1/\del) + C \log\log(e/\del), \del \Br) 
		+ C \log(e/\del)^{-1} } \\
	& = & (\kap W \log(1/\del))^{1/\kap} + O \Bl( \log\log(e/\del) \log(e\del)^{1/\kap - 1} \Br) \\
	& = & (\kap W \log(1/\del))^{1/\kap} + o(1)	\quad \mbox{as } \del \searrow 0 .
\eea

\paragraph{\sc Proof}
Due to the factorization lemma, the conditional probability and expectation factorize under the above conditions, i.e. we may consider a sequence $(y_n)_{n\in\N}$ and work with the sequence of conditional laws $\LL(Z_n(.,Y_n)\arrowvert Y_n=y_n)$, but note that we do not assume equality of $\LL(Z_n(.;Y_n)\arrowvert Y_n=y_n)$ and $ \LL(Z_n(.;y_n))$ in general. 
The first part of the proof is a modification of the Chaining in D\"umbgen and Walther (2008, technical report) applied to the conditional distribution $\LL(Z_n(.,Y_n)\arrowvert Y_n=y_n)$ for $y_n\in\CC_n$. Here we need however to define their additive correction function $H_1$ ($=H_n(.;y_n)$ subsequently due to the dependence on $n$ and $Y_n$ in our setting) in a different way, taking into account the additional logarithmic terms in the bound of the covering numbers. 
Lining up with their arguments, a suitable choice for the correction function $H_n$ appears to be
\begin{align*}
G_n&\Bigg\{W\log\Big(\frac{1}{\sigma_n(t;y_n)}\Big)\, +\, (B+\alpha)\log u\big(\sigma_n(t;y_n)\big)
\ +\  (2+\alpha)\log\log\Big(\frac{e}{\sigma_n(t;y_n)}\Big),\ \sigma_n(t;y_n)\Bigg\}\\
&=\ 
G_n\Bigg\{W\log\Big(\frac{1}{\sigma_n(t;y_n)}\Big)
+\, \Big((B+\alpha)\gamma + (2+\alpha)\Big)\log\log\Big(\frac{e}{\sigma_n(t;y_n)}\Big),\ \sigma_n(t;y_n)\Bigg\}.
\end{align*}
This term is essential for our proof of efficiency. It is important that the constant $\alpha$ does not influence the leading term. 
Concerning the tightness in probability as stated in the second part of Theorem \ref{Levy}, notice that it does not follow by an immediate continuity argument because the metric (and the metric space) change with both, $Y_n$ and $n$, hence some additional uniformity is required. 
For $0\leq\delta<\delta'\leq 1$ let $U_n(\delta,\delta';Y_n)$ be defined by
\begin{align*}
&\sup_{\substack{ \sigma_n(t;Y_n)\in(\delta,\delta']\\t\in\TT_n}}\Bigg\{\frac{|Z_n(t;Y_n)|}{\sig_n(t;Y_n)}  -   G_n \biggl( W
        \log\big(1/\sig_n(t;Y_n)\big) +  C \log\log\big(e/\sigma_n(t;Y_n)\big),\, \sigma_n(t;Y_n) \biggr)\Bigg\}.
\end{align*}
First observe that for any fixed $K>0$,
\begin{equation}\label{eq: ende}
\P\Big(U_n(0,1;Y_n) >K\Big\arrowvert Y_n\Big)\ \leq\ \P\Big(U_n(0,\delta;Y_n)>K/2\Big\arrowvert Y_n\Big)\ +\ \P\Big(U_n(\delta,1;Y_n)>K/2\Big\arrowvert Y_n\Big).
\end{equation}
The first part of Theorem \ref{Levy} implies that the first term on the right-hand-side in (\ref{eq: ende}) is bounded by $1-\P\big(\AA_n(2\del,q,Q;Y_n)\big\arrowvert Y_n\big)+C\log(e/\delta)^{-1}$ for $K>2C\log(e/\delta)^{-1}$ whenever $Y_n\in\CC_n$. Concerning the second term in (\ref{eq: ende}), note that
$$
U_n(\delta,1;Y_n)\ \leq\ -\inf_{\delta'\in[\delta,1]}\ H_n(\delta';Y_n)\ +\ \frac{1}{\delta}\sup_{\substack{t\in\TT_{Y_n}:\\ \sigma_n(t;Y_n)\geq\delta}}\Big\arrowvert Z_n(t;Y_n)\Big\arrowvert.
$$
Then the conclusion follows if we establish that
$$
\lim_{K\rightarrow\infty}\limsup_{n\rightarrow\infty}\P\bigg(\sup_{t\in\TT_{Y_n}}\big\arrowvert Z_n(t;Y_n)\big\arrowvert\, >\, K;\, Y_n\in\CC_n\bigg\arrowvert Y_n\bigg)\ =\ 0\ \ a.s.
$$
For $\varepsilon>0$ and $y_n\in\CC_n$, let $t_1(y_n),...,t_{m(y_n)}(y_n)$ be a maximal subset of $\TT_{y_n}$ with $\rho_n(t_i,t_j;y_n)>\varepsilon$ for arbitrary different indices $i,j\in\{1,...,m(y_n)\}$. Note that $m(y_n)\leq A\varepsilon^{-B}\log(e/\varepsilon)^{\alpha}$ by assumption (ii). Then condition (i) implies that
\begin{equation}\label{eq: 9}
\lim_{K\rightarrow\infty}\limsup_{n\rightarrow\infty}\P\bigg(\sup_{i=1,...,m(y_n)}\big\arrowvert Z_n\big(t_i({y_n});Y_n\big)\big\arrowvert >K\bigg\arrowvert Y_n=y_n\bigg)\ =\ 0\ \ a.s.
\end{equation}
On the other hand, we have on the set $\AA_n(\varepsilon,q,Q;Y_n)$ the bound
\begin{align}
\sup_{t\in\TT_{Y_n}}\arrowvert Z_n(t;Y_n)\arrowvert\
&\leq \ Q\varepsilon\log (e/\varepsilon)^{q}\  +\  \sup_{i=1,...,m(Y_n)}\big\arrowvert Z_n(t_i({Y_n});Y_n)\big\arrowvert.\label{eq: 9'}
\end{align} 
With $\varepsilon$ tending to zero sufficiently slowly, (\ref{eq: 9}) and (\ref{eq: 9'}) show together with the stochastic equicontinuity condition 
$\lim_{\delta\searrow 0}\inf_n\P\Big(\AA_n(\del,q,Q;Y_n)\Big\arrowvert Y_n\Big)=1$ a.s. 
\begin{align*}
\lim_{K\rightarrow\infty}\limsup_{n\rightarrow\infty}\P\bigg(\sup_{t\in\TT_{y_n}}\big\arrowvert Z_n(t;Y_n)\big\arrowvert > K\bigg\arrowvert Y_n=y_n\bigg)\ =\ 0\ \ a.s.
\end{align*}
Since the assumption $\inf_n\sup_{t\in\TT_{Y_n}}\sigma_n(t;Y_n)>0$ a.s. guarantees 
$$
\lim_{K\rightarrow\infty}\, \sup_n\P\Big(U_n(Y_n)<-K\Big\arrowvert Y_n\Big)\ =\ 0\ \ a.s.,
$$
the tightness in $(\P^{Y_n})$-probability is proved.\hfill$\square$

\medskip
\subsection{Proofs of the main results}
\paragraph{\sc Proof of Theorem \ref{Thm: limit distribution}}
Let $\lambda_n:=m/n$. In view of the $T_{jkn}$'s, the behavior of the process
$$
\Bigg(\frac{\sqrt{\lambda_n(1-\lambda_n)}}{\sqrt{n}}\sum_{i=0}^k\psi\bigg(\frac{\Arrowvert X_j-X_j^i\Arrowvert_2}{\Arrowvert X_j-X_j^k\Arrowvert_2}\bigg)\big(\Lambda\circ\Pi\big)(X_j^i)\Bigg)_{1\leq j\leq n,\, 0< k\leq n-1}
$$ 
conditional on $\XX_n$ needs to be investigated, where
$\Lambda\circ\Pi\arrowvert\XX_n$ is uniformly distributed on the set
$$
\bigg\{\lambda:\XX_n\rightarrow\big\{1/\lambda_n,-1/(1-\lambda_n)\big\}:\sum_{x\in\XX_n}\lambda(x)=0\bigg\}.
$$
For notational convenience it seems useful to redefine the process on
the random index set 
$$ 
\hat{\TT}_n\ :=\ \Big\{ \big(X_j,\Arrowvert X_j-X_j^k\Arrowvert_2\big):\, 1\leq j\leq n,\ 0 < k\leq
n-1\Big\}
$$ 
via the map 
$
(j,k)\ \mapsto \big(X_{j}, \,\big\Arrowvert X_{j}-X_{j}^k\big\Arrowvert_2\big)
$
and extend it to a process $\big(Y_n(t,r)\big)_{(t,r)\in\TT}$ with $\TT:=\big\{(t,r):\, t\in[0,1]^d, \,
0<r\leq\max_{x\in[0,1]^d}\Arrowvert x-t\Arrowvert_2\big\}$ by the definition
$$
Y_n(t,r)\ :=\
{\sqrt{n}}\sqrt{\lambda_n(1-\lambda_n)}\int\psi\bigg(\frac{\Arrowvert t-x\Arrowvert_2}{r}\bigg)\Big(d\hat{\P}_n^{\Pi}(x)-d\hat{\Q}_n^{\Pi}(x)\Big),
$$
where $\hat{\P}_n^{\Pi}$ and $\hat{\Q}_n^{\Pi}$ denote the empirical measures based on the permutated variables $X_{\Pi(1)},..., X_{\Pi(m)}$ and $X_{\Pi(m+1)},..., X_{\Pi(n)}$, respectively.
Let
\begin{align*}
\hat{\gamma}_n(t,r)^2\ :&=\ \mathrm{Var}\Big(Y_n(t,r)\Big\arrowvert
\XX_n\Big)\\
&=\ \frac{n}{n-1}\int\bigg[\psi\bigg(\frac{\Arrowvert t-x\Arrowvert_2}{r}\bigg)-\int\psi\bigg(\frac{\Arrowvert t-z\Arrowvert_2}{r}\bigg)d\hat{\H}_n(z)\bigg]^2d\hat{\H}_n(x),
\end{align*}
with $\hat{\H}_n$ the empirical measure of the observations $X_1,...,X_n$.

\medskip
\noindent
In the sequel we make use of the results in the previous section twice - in order to prove the tightness and weak approximation in probability of the sequence of conditional test statistics and within the "loop" we use the chaining arguments again to establish a sufficiently tightened uniform stochastic bound for the covering numbers below.

\medskip
\noindent
{\sc I. (Subexponential increments and Bernstein type tail behavior)}
The inversion of the conditional Bernstein type exponential inequality in Proposition
\ref{Prop: Bernstein} shows that for any
$\eta>0$, 
$$
\P\biggl(\Big\arrowvert\frac{Y_n(t,r)}{\hat{\gamma}_n(t,r)}\Big\arrowvert>G_n\big(\eta,\hat{\gamma}_n(t,r)\big)\bigg\arrowvert
{\XX}_n\biggr)\ \leq\ 2\,\exp(-\eta),
$$
where
$$
G_n\big(\eta,\hat{\gamma}_n(t,r)\big)\ :=\ R_n\big(\hat{\gamma}_n(t,r)\big)\eta\ +\ \Bigl(\big(R_n\big(\hat{\gamma}_n(t,r)\big)\eta\big)^2+2\delta(m,n)^2\eta\Bigr)^{1/2}\ \
\ \ \ \ \ \ \ \ \ \ \ \
$$  
with 
$$
R_n(\tau)\ :=\ \delta(m,n)\frac{2\Arrowvert
  \psi\Arrowvert_{\sup}\sqrt{\lambda_n(1-\lambda_n)}}{3\,\min(\lambda_n, 1-\lambda_n)\sqrt{n}\,\tau}.
$$
Let the random pseudo-metric $\hat{\rho}_{n}$ on $\TT$ be defined by
\begin{align*}
\hat{\rho}_n\bigl((t,r),(t',r')\big)^2
:=&\
\mathrm{Var}\Bigl(Y_n(t,r)-Y_n(t',r')\Big\arrowvert {\XX}_n\Bigr)\\
=& \frac{n}{n-1}\bigg[\int\Big(\psi_{tr}(x)-\psi_{t'r'}(x)\Big)^2d\hat{\H}_n(x) \ -
\ \bigg(\int\Big(\psi_{tr}(x)-\psi_{t'r'}(x)\Big)d\hat{\H}_n(x)\bigg)^2\bigg],
\end{align*}
with
$
\psi_{tr}(x)\ :=\ \psi\big(\frac{\Arrowvert t-x\Arrowvert_2}{r}\big).
$
Then the application of the second exponential
inequality of Proposition \ref{Prop: Bernstein} implies for any fixed $(t,r),
(t',r')\in\TT$ that
\begin{align*}
\P\biggl(\big\arrowvert Y_n(t,r)-Y_n(t',r')\big\arrowvert >
\hat{\rho}_{n}\big((t,r),(t',r')\big)\,q\,\eta\bigg\arrowvert {\XX}_n\biggr)
\leq\ 2\,\exp(-\eta),
\end{align*}
where 
$$
q\ :=\ 2\biggl(1\ +\ \frac{9\lambda_n(1-\lambda_n)}{2\max(\lambda_n,1-\lambda_n)^2}\big(\log 2\big)^{-1}\biggr).
$$

\medskip
\noindent
{\sc II. (Random local covering numbers)}  We need a bound for the local random covering numbers
$
N\bigl( (u\delta)^{1/2},\big\{ (t,r)\in \hat{\TT}_n: \hat{\gamma}_n(t,r)^2\, \leq\, \delta\big\},\hat{\rho}_n\bigr).
$
This is the most involved part of the proof. In contrast to previous work we aim at a uniform stochastic bound. In order to establish a sufficiently sharp upper bound, the following two claims are established:

\medskip
\noindent
(i) Let 
$$
\hat{\rho}_{2,n}\big((t,r),(t',r')\big)^2\ :=\ \int\Big(\psi_{tr}(x)-\psi_{t'r'}(x)\Big)^2d\hat{\H}_n(x)
$$
and define $d_n$ for arbitrary different points in $\hat{\TT}_n$ via
$$
{d_n}^2 := 
\max\big[\E\,{\hat{\rho}_{2,n}}^2,4/n\big]\bigg(1+C\log\Big(4\,e\Big/\max\big[\E\,{\hat{\rho}_{2,n}}^2,4/n\big]\Big)\bigg),
$$ 
with $C$ a positive constant to be chosen later.
Note that the map $x\mapsto x\sqrt{1+2C\log(\sqrt{e}/x)}$ is subadditive for $x\in(0,1]$, hence $d_n$ defines a metric. Furthermore let ${\gamma_n}^2:=\E\,\hat{\gamma}_{2,n}^2-\big(\E\hat{\gamma}_{1,n}\big)^2$, where
$$
\hat{\gamma}_{1,n}(t,r)^2\ :=\ \bigg(\int\psi_{tr}(x)d\hat{\H}_n(x)\bigg)^2\ \ \text{and}\ \  \hat{\gamma}_{2,n}(t,r)^2\ :=\ \int\psi_{tr}(x)^2d\hat{\H}_n(x).
$$
Then there exist a constant $C'>0$ and  a sequence $(\CC_n)_{n\in\N}$ of measurable sets with $\P_n^{\otimes m}\otimes\Q_n^{\otimes (n-m)}(\CC_n)\rightarrow 1$, such that for any  $\delta>0$, $u\in (0,1]$ with $u\delta\geq 4/n$ and any realization  $(X_1,...,X_n)\in\CC_n$ 
\begin{align*}
N\Bigl( (u\delta)^{1/2},\Big\{ (t,r)\in \hat{\TT}_n:\, &\hat{\gamma}_n(t,r)^2\, \leq\, \delta\Big\},\hat{\rho}_n\Bigr)\\ 
&\leq\  \ 
N\Bigl((u\delta)^{1/2},\Big\{(t,r)\in\hat{\mathcal{T}}_n:\, \gamma_{2,n}(t,r)^2\, \leq\,
C'\delta\log(e/\delta)^4\Big\},d_n\Bigr),  
\end{align*}
if $\psi$ is not rectangular. In case of the rectangular kernel, the set 
$$
\Big\{(t,r)\in\hat{\mathcal{T}}_n: \gamma_{2,n}(t,r)^2\, \leq\,
C'\delta\log(e/\delta)^4\Big\}
$$ in the covering number has to be replaced by
$$
\Big\{(t,r)\in\hat{\mathcal{T}}_n:\gamma_{2,n}(t,r)^2\ \leq\ C'\delta\log\big(e\big/\delta\big)^4\Big\}\cup\Big\{(t,r)\in\hat{\mathcal{T}}_n:\gamma_{2,n}(t,r)^2\ \geq\ 1-C'\delta\log\big(e\big/\delta\big)^4\Big\}.
$$
\noindent
(ii) There exists a constant $A>0$, independent of $u,\delta$ and $n$, such that whenever $u\delta\geq 4/n$, the upper bound given in (i) is again bounded from above by 
$A u^{-(d+1)}\delta^{-1}\log\big(e/(u\delta)\big)^{5(d+1)}. 
$
Moreover, the latter bound remains valid with $\TT$ in place of $\hat{\TT}_n$.

\medskip
\noindent
Note that we cannot rely our bound directly on uniform covering numbers and Vapnik-Cervonenkis (VC) theory as the envelope $I\{X\in\XX_n\}$ only allows for a bound of order $u^{-2}\delta^{-2}$, which would result in the loss of efficiency of the procedure, and a pre-partitioning of $\hat{\TT}_n$ as used in the proof of (ii) seems to be rather involved.

\medskip
\noindent 
{\it Proof of (i):~~} We first derive a uniform stochastic bound for the random metric $\hat{\rho}_{2,n}$. Recall that every function $\psi$ of bounded total variation is representable as a difference of isotonic functions $\psi^{(1)}$ and $\psi^{(2)}$. With the definition of the subgraphs 
$$
\text{sgr}\big(\psi^{(i)}_{tr}\big)\ :=\ \Big\{ (x,y)\in[0,1]^d\times \R:\, y\leq \psi^{(i)}_{tr}(x)\Big\},\ \ i=1,2,
$$
the set $\big\{\text{sgr}\big(\psi_{tr}^{(i)}\big): (t,r)\in\TT\big\}$ has a  VC-dimension bounded by $d+3$ (van der Vaart and Wellner 1996) with envelope $TV(\psi)$.  
Consequently, the uniform covering numbers $N(\varepsilon,\FF)$ with
$$
\FF\ :=\ \Big\{\big(\psi_{tr}-\psi_{t'r'}\big)^2:\, (t,r),(t',r')\in\TT\Big\}
$$ 
is bounded by $C\varepsilon^{-\alpha}$ for some real-valued $\alpha>0$ and some constant $C>0$. The boundedness of $\psi$ shows that $\FF$ is uniform Glivenko-Cantelli in particular (see Dudley, Gin\'e and Zinn 1991, for instance). As an immediate consequence,
\begin{equation}
\lim_{n\rightarrow\infty}\P\bigg(\Big\Arrowvert \hat{\rho}_{2,n}\big((t,r),(t',r')\big)^2-\E\hat{\rho}_{2,n}\big((t,r),(t',r')\big)^2\Big\Arrowvert_{\TT\times\TT}>\delta\bigg)\ =\ 0, \label{eq: unif. conv.}
\end{equation}
for any $\delta>0$. However such a bound is not sufficient for our purposes. Because of $\Arrowvert \psi\Arrowvert_{\sup}\leq 1$, the squared random metric 
 ${\hat{\rho}_{2,n}}^2$ is $1/n$ times the sum of $n$ independent random variables with absolute values $\leq 4$, hence
\begin{align*}
\text{Var}\Big({\hat{\rho}_{2,n}\big((t,r),(t',r')\big)}^2\Big)\ \leq\ \frac{4}{n}\,\E\Big(\hat{\rho}_{2,n}\big((t,r),(t',r')\big)^2\Big)\ 
&\leq\ \max\bigg\{\frac{4}{n},\E\Big(\hat{\rho}_{2,n}\big((t,r),(t',r')\big)^2\Big)\bigg\}^2.
\end{align*}
Now the application of Bernstein's exponential inequality (see Shorack and Wellner 1986) entails

\vspace{-3.7mm}
\begin{align*}
\P\bigg(\bigg\arrowvert\frac{\hat{\rho}_{2,n}\big((t,r),(t',r')\big)^2-\E\hat{\rho}_{2,n}\big((t,r),(t',r')\big)^2}{\max[4/n,\E\hat{\rho}_{2,n}\big((t,r),(t',r')\big)^2]}\bigg\arrowvert\ >\ \eta\bigg)\ 
&\leq\ 2\exp\Bigg(-\frac{\eta^2/2}{1+\eta/3}\Bigg)\\
&\leq\ 2\exp\bigg(-\frac{3}{2}\eta+\frac{9}{2}\bigg)
\end{align*}
for arbitrary points $(t,r),(t',r')\in\TT$. I.e. $\hat{\rho}_{2,n}^2-\E\hat{\rho}_{2,n}^2$, standardized by $
\max\big\{4/n,\E\hat{\rho}_{2,n}^2\big\}
$, 
has (uniformly) subexponential tails. Analogously, the process $\hat{\rho}_{2,n}^2-\E\hat{\rho}_{2,n}^2$ has subexponential increments with respect to the metric $\tilde{D}_n$ given by $$
\tilde{D}_n\big(a,b\big)\ :=\ \max\Big[1/n,\E\Big(\hat{\rho}_{2,n}^2(a)-\hat{\rho}_{2,n}^2(b)\Big)^2\Big]I\big\{a\not=b\big\},\ \ a,b\in\TT\times\TT.
$$ 
Note that $
\max[4/n,\E\hat{\rho}_{2,n}^2]$ is Lipschitz continuous with respect to $\tilde{D}_n$.
Theorem \ref{Chaining} shows that the above ingredients imply that $\lim_{\delta\searrow 0}\inf_n\P\big(\AA_n(\del,1,Q;\XX_n)\big\arrowvert \XX_n\big)=1$ for some adequately chosen $Q>0$, where we use  the definition of $\AA_n$ from Theorem \ref{Levy} with $Y_n=\XX_n$ and $Z_n=\hat{\rho}_{2,n}^2-\E\hat{\rho}_{2,n}^2$. Now we may apply the latter
to conclude that there exists some universal constant $C>0$ such that the probability of the event
\begin{align}\label{eq: unif. conv. 2}
\bigg\{\Big\arrowvert \hat{\rho}_{2,n}\big(&(t,r),(t',r')\big)^2-\E\hat{\rho}_{2,n}\big((t,r),(t',r')\big)^2\Big\arrowvert >\\ 
&\ \ \  \ \ \ C\,\max\big[4/n,\E\hat{\rho}_{2,n}\big((t,r),(t',r')\big)^2\big]\log\Big(4\,e\Big/\max\big[4/n,\E\hat{\rho}_{2,n}\big((t,r),(t',r')\big)^2\big]\Big)\nonumber\\
& \text{~~for some~}(t,r),(t',r')\ \text{with~}\E\hat{\rho}_{2,n}\big((t,r),(t',r')\big)^2\leq\delta'\bigg\}\nonumber
\end{align}
is bounded by some function $\varepsilon(\delta')$ independent of $n$ with $\lim_{\delta'\searrow 0}\varepsilon(\delta')=0$. Since the probability in (\ref{eq: unif. conv.}) is antitonic in $\delta$ for any fixed $n$ with limes $0$ as $n\rightarrow\infty$ for any fixed $\delta$, there exists a sequence $\delta_n\searrow 0$ along which the result of (\ref{eq: unif. conv.}) still holds true. Thus,  combining (\ref{eq: unif. conv.}) and (\ref{eq: unif. conv. 2}) for a sequence $\delta'=\delta_n'\searrow 0$ sufficiently slowly implies the existence of a sequence of sets $(\AA_n)_{n\in\N}$ with $\P^{\otimes m}\otimes\Q^{\otimes (n-m)}(\AA_n)\rightarrow 1$ such that 
$$
\hat{\rho}_{2,n}\ \leq\ \max\big[4/n,\E\hat{\rho}_{2,n}^2\big]^{1/2}\Big(1+C\log\big(4\,e\big/\max\big[4/n,\E\hat{\rho}_{2,n}^2]\big)\Big)^{1/2}
\ \ \text{whenever $\underbar{X}\in\AA_n$.}
$$
The treatment of the random set
$$
\hat{\BB}_{\delta}\ :=\ \Big\{ (t,r)\in \hat{\TT}_n: \hat{\gamma}_n(t,r)^2\, \leq\, \delta\Big\}
$$
is similar in spirit but more involved because the random quantity $\hat{\gamma}_n^2$ is not representable as a sum of independent variables. However we can use the decomposition $[(n-1)/n] {\hat{\gamma}_n}^2={\hat{\gamma}_{2,n}}^2-{\hat{\gamma}_{1,n}}^2$. Before deriving a stochastic bound, we notice the following: If $\psi$ describes the rectangular kernel, we have $\hat{\gamma}_{2,n}^2=\hat{\gamma}_{1,n}$, i.e. 
$$
\hat{\gamma}_{2,n}^2-\hat{\gamma}_{1,n}^2\ =\ \hat{\gamma}_{2,n}^2\big(1-\hat{\gamma}_{2,n}^2\big).
$$
In this case, the random set $\hat{\BB}_{\delta}$ is consequently contained in the union
\begin{equation}\label{eq: union}
\Big\{\hat{\gamma}_{2,n}^2\leq 2\delta\Big\} \cup \Big\{\hat{\gamma}_{2,n}^2\geq 1-2\delta\Big\}.\end{equation}
Consider the general case. Using that
\begin{equation}\label{eq: varianz 1}
\text{Var}\Big({\hat{\gamma}_{1,n}(t,r)}\Big)\ =\ \frac{1}{n^2}\sum_{i=1}^n\Big(\E\psi_{tr}(X_i)^2-\big(\E\psi_{tr}(X_i)\big)^2\Big)\
\leq\ \frac{1}{n}\E\Big(\hat{\gamma}_{2,n}(t,r)^2\Big)
\end{equation}
and
\begin{equation}\label{eq: varianz 3}
\text{Var}\Big({\hat{\gamma}_{2,n}(t,r)}^2\Big)\ =\ \frac{1}{n^2}\sum_{i=1}^n\Big(\E\psi_{tr}(X_i)^4-\big(\E\psi_{tr}(X_i)^2\big)^2\Big)\
\leq\ \frac{1}{n}\E\Big(\hat{\gamma}_{2,n}(t,r)^2\Big),
\end{equation}
we may apply the above chain of arguments for $\hat{\rho}_{2,n}^2$ to $\hat{\gamma}_{1,n}$ and $\hat{\gamma}_{2,n}^2$ together with the upper bounds in (\ref{eq: varianz 1}) and (\ref{eq: varianz 3}) for the standardization respectively and 
obtain the existence of a constant $C_1>0$ such that
\begin{align*}
\gamma_{1,n}\ -\ &\frac{C_1\max\big[1/n,\gamma_{2,n}^2\big]^{1/2}}{\sqrt{n}}\log\Big(e\sqrt{n}\Big/\max\big[1/n,\gamma_{2,n}^2\big]^{1/2}\Big)\\ 
&\leq\ \hat{\gamma}_{1,n}\ \leq\ \gamma_{1,n}\ +\ \frac{C_1\max\big[1/n,\gamma_{2,n}^2\big]^{1/2}}{\sqrt{n}}\log\Big(e\sqrt{n}\Big/\max\big[1/n,\gamma_{2,n}^2\big]^{1/2}\Big)
\end{align*}
whenever $\underbar{X}\in\DD_n$ for some sequence $(\DD_n)_{n\in\N}$ with asymptotic probability $1$, uniformly evaluated at $(t,r)\in\hat{\TT}_n$. Note that $\hat{\gamma}_{1,n}\geq 1/n$, $\hat{\gamma}_{2,n}^2\geq 1/n$ for all $(t,r)\in\hat{\TT}_n$. 
The same holds true with a constant $C_2>0$ and a sequence $(\DD_n')_{n\in\N}$ with asymptotic probability $1$ and $\hat{\gamma}_{1,n}$ and $\gamma_{1,n}$ replaced by $\hat{\gamma}_{2,n}^2$ and $\gamma_{2,n}^2$. Using the lower bound for $\hat{\gamma}_{2,n}^2$ and the upper bound for $\hat{\gamma}_{1,n}$, a bit of  algebra yields 
\begin{align*}
\hat{\BB}_{\delta}\ \subset\ \bigg\{\gamma_{2,n}^2-\gamma_{1,n}^2\
&\leq\ \delta\ +\ \max\big[1/n,\gamma_{2,n}^2\big]^{1/2}\frac{K}{\sqrt{n}}\log\Big(e\sqrt{n}\Big/\max\big[1/n,\gamma_{2,n}^2\big]^{1/2}\Big)^2\bigg\}
\end{align*}
whenever $\underbar{X}\in \DD_n\cap\DD_n',\ \delta\geq 1/n$. Here and from now on, $K$ denotes some universal constant, not dependent on $n$ and $(t,r)$. Its value may be different in different expressions. Now we first consider the case
$$
\sup_{n\in\N}\sup_{(t,r)\in\TT}\Big(\gamma_{1,n}^2\Big/\gamma_{2,n}^2\Big)\ \leq\ C'\ <\  1.
$$
Then the above condition shows that
\begin{align*}
\gamma_{2,n}^2(1-C')\ &\leq\ \delta\ + \max\big[1/n,\gamma_{2,n}^2\big]^{1/2}\frac{K}{\sqrt{n}}\log\Big(e\sqrt{n}\Big/\max\big[1/n,\gamma_{2,n}^2\big]^{1/2}\Big)^2\\
&\leq\ 2\max\bigg\{\delta,\ \max\big[1/n,\gamma_{2,n}^2\big]^{1/2}\frac{K}{\sqrt{n}}\log\Big(e\sqrt{n}\Big/\max\big[1/n,\gamma_{2,n}^2\big]^{1/2}\Big)^2\bigg\},
\end{align*}
which entails that $\gamma_{2,n}^2\leq\ K\,\delta\log\big(e\big/\delta\big)^4$ for $\delta\geq 1/n$ by  the isotonicity of $x\mapsto x\log(e/x)^4$ on $(0,1]$. 
On the other hand, the case
\begin{equation}\label{eq: AA}
\sup_{n\in\N}\sup_{(t,r)\in\TT}\Big(\gamma_{1,n}^2\Big/\gamma_{2,n}^2\Big)\ =\  1
\end{equation}
implies already that $\psi$ is equal to the rectangular kernel: If the sup is attained it is obvious. The equicontinuity of $(h_n)_{n\in\N}$ and its uniformly bounded $L_1$-norm $\Arrowvert h_n\Arrowvert_1=1$ imply its uniform boundedness, hence relative compactness in the  topology of uniform convergence by the Arzel$\grave{\text{a}}$-Ascoli-Theorem. There therefore exists  at least  a uniformly convergent subsequence $(h_{m(n)})$ with (uniformly) continuous limit, say $h$,  along this result holds true as well, because
$
\max_{(t,r)\in\TT}\big(\gamma_{1,n}^2\big/\gamma_{2,n}^2\big)
$
depends continuously on the mixed density. This however implies that $\psi$ describes the rectangular kernel, because the uniform limit $h$ of that subsequence is bounded away from zero. Hence in case of (\ref{eq: AA}), we consequently  obtain by (\ref{eq: union})
\vspace{-2mm}
$$
\hat{\BB}_{\delta}\ \subset\ \Big\{\gamma_{2,n}^2\ \leq\ K\delta\log\big(e\big/\delta\big)^4\Big\}\cup\Big\{\gamma_{2,n}^2\ \geq\ 1-K\delta\log\big(e\big/\delta\big)^4\Big\}  \ \text{whenever~}\underbar{X}\in\DD_n\cap\DD_n',\ \delta\geq 1/n.
$$

\noindent
{\it Proof of (ii):~~} 
Since $\psi$ is of bounded total
variation, there exists some finite measure $\mu$ such that for any $0\leq
z_1<z_2\leq 1$, $\arrowvert\psi(z_1)-\psi(z_2)\arrowvert \leq
\mu[z_1,z_2]$. With
$$
M_x(t,t',r,r')\ :=\ \bigg[0,\frac{\Arrowvert t-x\Arrowvert_2}{r}\bigg]\Delta\bigg[0,\frac{\Arrowvert t'-x\Arrowvert_2}{r'}\bigg]
$$

\vspace{-2mm}
\noindent
we obtain
\vspace{-2mm}
\begin{align}
\E\hat{\rho}_{2,n}\big((t,r),(t',r')\big)^2\ &\leq\ \int\bigl(\psi_{tr}(x)-\psi_{t'r'}(x)\bigr)^2d{\H_n}(x)\nonumber\\
&\leq\ K\int\mu\bigl(M_x(t,t',r,r')\big)\,d\H_n(x)\nonumber\\
&=\ K\int\int I\big\{y\in M_x(t,t',r,r')\}\,d\H_n(x)d\mu(y)\ \ \ \ \  \text{(Fubini)}\nonumber\\
&\leq \ K\sup_{y\in[0,1]}\int I\big\{ y\in M_x(t,t',r,r')\big\}\,d\H_n(x).\label{eq: A}
\end{align}
Then $y\in
M_x(t,t',r,r')$ implies that 
$
x\in B_t\big(ry\big)\Delta B_{t'}\big(r'y\big)$. Since $h_n$ is uniformly
bounded from  above, we obtain that (\ref{eq: A}) is not greater than $ K\lambda\big(B_t(r)\Delta B_{t'}(r')\big)$. Consequently,  $d_n\leq K\,d$ if $d_n\geq 4/n$ with the metric $d$ defined below in (\ref{metric}), due to the isotonicity of $x\mapsto x(1+C\log(e/x))$ for $x\in (0,1]$, $C>0$. 
$\psi$ attains its maximum $1$ at $0$, hence there exists some $r^*>0$ such that $\psi(\Arrowvert x\Arrowvert_2)\geq 1/2$ whenever $\Arrowvert x\Arrowvert_2\leq r^*$.
 Using in addition the uniform boundedness of $h_n$ away from zero we obtain
$
\gamma_{2,n}(t,r)^2 \geq  K\cdot r^d\ \ (t,r)\in\TT.
$
We now start bounding the covering numbers
\vspace{-2mm}
$$
N\bigg((u\delta)^{1/2},\, \Big\{(t,r)\in{\mathcal{T}}: \gamma_{2,n}(t,r)^2\leq 
K\delta\log(e/\delta)^4\Big\},\, d\bigg),
$$
where the metric $d$ on $\TT\times\TT$ is pointwise defined by
\begin{equation}
d\big((t,r),(t',r')\big)^2\ :=\ \lambda\big(B_t(r)\Delta B_{t'}(r')\big)\bigg(1+C\log\Big[V\,e\Big/\lambda\big(B_t(r)\Delta B_{t'}(r')\big)\Big]\bigg)\label{metric}
\end{equation}
with $V=\lambda(B_0(\sqrt{d}))$ the volume of the $d$-dimensional Euclidean ball with radius $\sqrt{d}$. 
Again by the isotonicity of $x\mapsto x\log(e/x)$ for $x\in (0,1]$, the inequality  $\tilde{d}\big((t,r),(t',r')\big):=\lambda\big(B_t(r)\Delta B_{t'}(r')\big)^{1/2}\leq\varepsilon/\sqrt{\log(V\,e/\varepsilon^2)}$ implies that $d\big((t,r),(t',r')\big)$ is not greater than $(2C+1)^{1/2}\varepsilon$. Thus in order to finish claim (ii), it is sufficient to bound 
\begin{equation}\label{eq: covering number}
N\bigg(\Big(\frac{u\delta}{\log(e/(u\delta))}\Big)^{1/2},\, \Big\{(t,r)\in\mathcal{T}: r^d\leq
\delta\log(e/\delta)^4\Big\},\, \tilde{d}\bigg).
\end{equation}

\vspace{-1mm}
\noindent
First note that there exists a finite collection of at most $m\leq K/(\delta\log(e/\delta)^4)$ points $t_1,...,t_m$ such that the set
$
\Big\{(t,r)\in\mathcal{T}: r^d\leq
\delta\log(e/\delta)^4\Big\}
$
is contained in the union $\cup_{i=1}^m\AA_i$ with
\vspace{-2.5mm}
$$
\AA_i \ :=\ \bigg\{(t,r)\in\TT: B_t(r)\subset B_{t_i}\Big([K'\delta\log(e/\delta)^4]^{1/d}\Big)\bigg\}
$$

\vspace{-0.5mm}
\noindent
for some universal $K'>0$. The rotation and translation invariance of the Lebesgue measure leads to the rescaling invariance for the covering numbers
\vspace{-2mm}
\begin{equation}
N\Big(\varepsilon^{1/2},\, \big\{(t,r):\, B_t(r)\subset B_{0}(R)\big\},\, \tilde{d}\,\Big)\ =\ N\Big((\varepsilon/R^d)^{1/2},\, \big\{(t,r):\, B_t(r)\subset B_{0}(1)\big\},\, \tilde{d}\,\Big).\label{cn}
\end{equation}

\vspace{-2mm}
\noindent
But a minimal $\tilde{d}$-$~(\varepsilon/R^d)^{1/2}$-net of the set $\big\{(t,r)\in\TT: B_t(r)\subset B_0(1), r=r'\big\}$ for some {\it fixed}  $r'> \varepsilon^{1/d}/R$ contains not more than $M=K[R^d/\varepsilon]^d$ elements $(t_1,r'),...,(t_M,r')$ with $K$ uniformly in $r'\in (\varepsilon^{1/d}/R,\sqrt{d}]$, noticing that $\lambda(B_t(r)\Delta B_{t'}(r))\leq K\Arrowvert t-t'\Arrowvert_2r^{d-1}$ and $r\leq \sqrt{d}$.  Now fix a $K(\varepsilon/R^d)$-net $t_1,...,t_M$ with respect to $\Arrowvert .\Arrowvert_2$ and observe that $\lambda(B_t(r)\Delta B_t(r'))\leq Kr^{d-1}(r-r')$ for $r>r'$, $r\leq\sqrt{d}$, which shows that the quantity (\ref{cn}) is bounded by $K(R^d/\varepsilon)^{d+1}$ (with $K$ uniformly in $\varepsilon$ and $R$). Correspondingly, this holds true for $N\big((u\delta/\log[e/(u\delta)])^{1/2},\AA_i, \tilde{d}\big)$, hence the covering number (\ref{eq: covering number}) is bounded by
 $A\delta^{-1}u^{-(d+1)}\log(e/u\delta)^{5(d+1)}$ for some universal constant $A>0$. An analogous bound holds for $\hat{\TT}_n$ in place of $\TT$ (and $u\delta\geq 4/n$): If $(t_1,r_1),...,(t_k,r_k)$ denotes an $\varepsilon$-net with respect to $d$ in $B\subset \TT$, we may define a $2\,\varepsilon$-net $(\hat{t}_1,\hat{r}_1),...,(\hat{t}_k,\hat{r}_k)$ in $\hat{\TT}_n\cap B$ via the definition $(\hat{t}_i,\hat{r}_i):=\argmin_{(t,r)\in\hat{\TT}_n\cap B}d\big((t,r),(t_i,r_i)\big)$. The corresponding covering numbers in case of the rectangular kernel for the sets $\big\{\gamma_{2,n}^2\, \geq\, 1-K\delta\log\big(e\big/\delta\big)^4\big\}$ can be treated with similar arguments, 
which concludes the proof of (ii).

\medskip
In order to line up with the requirements of Theorem \ref{Levy}, let us remark that the proof of that chaining requires only the special choice $u=u(\delta)=\log(e/\delta)^{\gamma}$ for some exponent $\gamma<0$, which entails that $\delta\leq n^{-1}(\log n)^{\alpha}$ for some $\alpha>0$ in case $u\delta\leq 4/n$. But for any $\alpha'>0$, $\sharp\big\{(t,r)\in\hat{\TT}_n:r^d\leq Kn^{-1}(\log n)^{\alpha'}\big\}= \sum_{i=1}^n\sharp\big\{(X_i,r)\in\hat{\TT}_n: r^d\leq Kn^{-1}(\log n)^{\alpha'}\big\}$, and with the same arguments as used in (i) we obtain for $r_n^d=n^{-1}(\log n)^{\alpha'}$ that the inequality $\hat{\H}_n\big(B_t(r_n))\leq K\lambda(B_t(r_n))\log n$ holds, uniformly in $t\in[0,1]^d$,  with asymptotic probability $1$, which entails $\sharp\big\{(t,r)\in\hat{\TT}_n:r^d\leq Kn^{-1}(\log n)^{\alpha'}\big\}=O_p\big(n(\log n)^{\alpha''}\big)$ for some $\alpha''>0$.

\medskip
\noindent
{\sc III. (Tightness and weak approximation in probability)} 
As a consequence of the above exponential inequalities in step I and the bound for the uniform
 covering numbers $N(\delta,{\TT})$, Theorem \ref{Chaining} shows 
\begin{equation}\label{eq: ec}
\lim_{\delta\searrow
  0}\limsup_{n\rightarrow\infty}\P\Biggl(\,\sup_{\hat{\rho}_n((t,r),(t',r'))\leq\delta}\frac{\arrowvert Y_n(t,r)-Y_n(t',r')\arrowvert}{\hat{\rho}_n((t,r),(t',r'))\log\big(e\big/\hat{\rho}_n((t,r),(t',r'))\big)}\ >\ \varepsilon\Bigg\arrowvert \XX_n \Biggr)\ =\ 0,
\end{equation}
where the $\sup$ within the brackets is even running over elements of ${\TT}\times{\TT}$. Now the application of Theorem \ref{Levy} entails that $\LL\big(T_n\circ\Pi\big\arrowvert\XX_n\big)$ is tight in $\big(\P_n^{\otimes m}\otimes\Q_n^{\otimes (n-m)}\big)$-probability. What remains being proved is the weak approximation.
Starting from (\ref{eq: ec}), the uniform convergence (\ref{eq: unif. conv.}) 
implies in particular the asymptotic stochastic equicontinuity 
$$
\lim_{\delta\searrow
  0}\limsup_{n\rightarrow\infty}\E_{(p_n,q_n,\lambda_n)}^*\P^*\bigg(\sup_{{\rho}_n((t,r),(t',r'))\leq\delta}\big\arrowvert Y_n(t,r)-Y_n(t',r')\big\arrowvert>\varepsilon\bigg\arrowvert\XX_n\bigg)\ =\ 0\ \ \text{for all $\varepsilon>0$.}
$$
Since to any subsequence of the metric $\rho_n$ there exists some uniformly convergent subsubsequence as a consequence of the relative compactness of $(h_n)_{n\in\N}$ in the uniform topology, it suffices (via proof of  contradiction) for the weak approximation in probability 
$$
d_w\Big\{\LL\Big(\big(Y_n(t,r))_{(t,r)\in\TT}\Big\arrowvert\XX_n\Big),\
\LL\Big(\big(Z_n(t,r)\big)_{(t,r)\in\TT}\Big)\Big\}\longrightarrow_{\P_n^{\otimes m}\otimes\Q_n^{\otimes (n-m)}}\ 0
$$
to establish the convergence of finite dimensional
distributions. Here, $d_w$ is defined via the outer expectations $\E^*$. 
For let $\big\{(t_{1},r_{1}),...,(t_{k},r_{k})\big\}$ be a collection of points from
${\TT}$. 
Denote furthermore 
$
 a_{rt}(X_i):=n^{-1/2}\sqrt{\lambda_n(1-\lambda_n)}\,\psi_{tr}(X_i).
$
Then
$$
\Big(Y_n(t,r)\Big)_{(t,r)\in\TT}\ =\
\bigg(\sum_{i=1}^na_{rt}(X_i)\Lambda(t^i)\bigg)_{(t,r)\in\TT},
$$
with $t^i$ the $i$'th nearest-neighbor of $t$ within $\XX_n$. 
Let $\big(Z_n(t,r)\big)_{(t,r)\in\TT}$ be pointwise be defined by
$
Z_n(t,r) := \sqrt{\lambda_n(1-\lambda_n)}\int\phi_{rt}^{(n)}(x)\,dW(x).
$
Using that $2\,\cov(X_1,X_2)$ equals $\mathrm{Var}(X_1)+\mathrm{Var}(X_2)-\mathrm{Var}(X_1-X_2)$ for two random variables $X_1$ and $X_2$, one finds that $[(n-1)/n]\mathrm{cov}\,\big(Y_n(t,r),Y_n(t',r')\big\arrowvert\XX_n\big)$ equals
\begin{align}
 -\frac{1}{2}\int&\Big(\psi_{tr}(x)-\psi_{t'r'}(x)\Big)^2d\hat{\H}_n(x)+\frac{1}{2}\bigg(\int\Big(\psi_{tr}(x)-\psi_{t'r'}(x)\Big)d\hat{\H}_n(x)\bigg)^2\,  +\, \frac{1}{2}\int\psi_{tr}(x)^2d\hat{\H}_n(x)\label{covarianz}\\
    &-\, \frac{1}{2}\bigg(\int\psi_{tr}(x)d\hat{\H}_n(x)\bigg)^2 +\, \frac{1}{2}\int\psi_{t'r'}(x)^2d\hat{\H}_n(x)\, -\, \frac{1}{2}\bigg(\int\psi_{t'r'}(x)d\hat{\H}_n(x) \bigg)^2.\nonumber
\end{align}
Replacing the empirical measure $\hat{\H}_n$ by its expectation $\H_n$, the above six expressions in (\ref{covarianz}) coincide with the covariance $\cov\big(Z_n(t,r),Z_n(t',r')\big)$ of the limiting process $Z_n$. Define $\bar{a}_{r_jt_j}^{(n)}:=n^{-1}\sum_{i=1}^na_{r_jt_j}^{(n)}(X_i)$, $j=1,...,k$. Since
$$
\sum_{j=1}^k\frac{\max_i(a_{r_jt_j}^{(n)}(X_i)-\bar{a}_{r_jt_j}^{(n)})^2}{\sum_{i=1}^n(a_{r_jt_j}^{(n)}(X_i)-\bar{a}_{r_jt_j}^{(n)})^2}\
\longrightarrow_{\P_n^{\otimes m}\otimes \Q_n^{\otimes (n-m)}}\ 0\ \ \ \ (n\rightarrow\infty)
$$
and
$
\big\arrowvert\mathrm{cov}\big(Y_n(t,r),Y_n(t',r')\big\arrowvert\XX_n\big)-\mathrm{cov}\big(Z_n(t,r),Z_n(t',r')\big)\big\arrowvert\longrightarrow_{\P_n^{\otimes m}\otimes \Q_n^{\otimes (n-m)}}
          0
$
by an application of the weak law of large numbers for triangular arrays to each of the expressions in (\ref{covarianz}) separately, H\'ajek's Central Limit Theorem for permutation
statistics extended for the multivariate setting yields the desired weak convergence in probability of the finite dimensional distributions.
For  notational convenience, define
$$
T_n^{\Pi}(\delta,\delta')\ :=\ \sup_{\substack{(j,k):\\ \delta<\gamma_n(j,k)\leq\delta'}}\Big\{\big\arrowvert T_{jkn}\circ\Pi\big\arrowvert\ -\
C_{jkn}\Big\}
$$
and
$$
S_n(\delta,\delta')\ :=\ \sup_{\substack{(t,r):\\ \delta<\gamma_n(t,r)\leq\delta'}}\left\{ \frac{\big\arrowvert\int
  {\phi}_{rt}^{(n)}(x)\,dW(x)\big\arrowvert}{\gamma_n(t,r)}\
  -\ \sqrt{2\log\bigl(1/{\gamma_n(t,r)}^2\bigr)} \,\right\}.
$$
Since $\sup_{t\in\TT\setminus\hat{\TT}_n}d(t,\hat{\TT}_n)\rightarrow_{\P^{\otimes m}\otimes\Q^{\otimes(n-m)}}\, 0$ and 
$
\sup_{(j,k):\,\gamma_n(j,k)\geq \delta}\big\arrowvert
C_{jkn}-(2\,\Gamma_{jkn})^{1/2}\big\arrowvert  \rightarrow_{\P_n^{\otimes m}\otimes \Q_n^{\otimes (n-m)}} 0$ as $n\rightarrow\infty$,
it follows from the above established results that
$$
d_w\Big(\LL\big(T_n^{\Pi}(\delta,1)\big\arrowvert\XX_n\big),\
\LL\big(S_n(\delta,1)\big)\Big)\ \longrightarrow_{\P_n^{\otimes m}\otimes \Q_n^{\otimes (n-m)}} \ 0
$$\nopagebreak
for
any fixed $\delta\in (0,1]$. An
application of Theorem \ref{Levy} as well as its
subsequent Remark imply that
$$
\lim_{\delta\searrow
  0}\limsup_{n\rightarrow\infty}\,\E\,\P\big(T_n^{\Pi}(0,\delta)\geq\varepsilon\big\arrowvert\XX_n\big)\ =\
0 \ \ \textrm{and}\ \ \ \lim_{\delta\searrow
  0}\limsup_{n\rightarrow\infty}\P\big(S_n(0,\delta)\geq\varepsilon\big)\ =\ 0
$$
for any $\varepsilon >0$. Thus, because obviously 
$
\lim_{\delta\searrow 0}\limsup_{n\rightarrow\infty}\P\big(S_n(\delta,1)\leq -\varepsilon\big)\ =\ 0,
$
we obtain
$$
d_w\Big(\LL\big(T_n^{\Pi}(0,1)\big\arrowvert\XX_n\big),\
\LL\big(S_n(0,1)\big)\Big)\ \longrightarrow_{\P_n^{\otimes m}\otimes \Q_n^{\otimes (n-m)}} \ 0.\eqno{\Box}
$$

\medskip
\paragraph{\sc Proof of Theorem \ref{Thm: lower bound}}
Let $\CC$
be some compact rectangle of $J$. Fix $\beta >0$. For any integer $k>1$ let  
$\CC_{n,k}\subset \CC$ be some maximal subset of points such that $\Arrowvert
x-y\Arrowvert_2\geq 2k\delta_n$ and $B_{x}(k\delta_n)\subset \CC$ for
arbitrary different points $x,y\in\CC_{n,k}$. Then $\sharp\CC_{n,k}\sim (k\delta_n)^{-d}$. Now let $\phi_{x,n}$ be the solution of the subsequent optimization problem:

\vspace{1.2mm}
\noindent
($*$)~~Minimize $\Arrowvert g\Arrowvert_2$ under the constraints
$$
g\in\HH_d(\beta,L;\R^d),\ \ \supp(g)\subseteq B_{x}(k\delta_n),\ \ g(x) = L\delta_n^{\beta},\ \ \int g(z)\sqrt{h_n(z)}dz=0.
$$
These constraints define a closed and convex set in $L_2\big([0,1]^d\big)$ which is non-empty for $k$ sufficiently large (and uniformly in $n$ due to the equicontinuity of $(h_n)$ and the rescaling property, see subsequently to (\ref{eq: ORR}) below). Consequently in the latter case, the argmin $\phi_{x,n}$ exists and is unique. The resulting density candidates 
$$
p_{x,n}\ =\ h_n\cdot\Big(1+\big(1-(m/n)\big)\phi_{x,n}\big/\sqrt{h_n}\Big)\ \  \text{and}\ \  q_{x,n}\ =\ h_n\cdot\Big(1-(m/n)\phi_{x,n}\big/\sqrt{h_n}\Big)
$$ 
are non-negative and thus contained in $\FF_{h_n}^{(m,n)}$ as soon as additionally 
\begin{equation}
-\,\frac{\sqrt{h_n(.)}}{1-m/n}\ \leq\ \phi_{x,n}(.)\ \leq\ \frac{\sqrt{h_n(.)}}{m/n}\ \ \ \text{for all~}x\in\CC_n.\nonumber
\end{equation}
This is guaranteed for sufficiently large $n$ when sequence $(\delta_n)_{n\in\N}$ tends to zero.  
For
any statistical level-$\alpha$-test $\psi=\psi(\beta,L,h_n):\R^{d\times n}\rightarrow [0,1]$ for testing the hypothesis "$\phi=0$" it holds true that
\begin{align}
\min_{x\in\CC_{n,k}}\E_{(m,n,p_{x,n},q_{x,n})}\psi -\ \alpha\ &\leq\ \min_{x\in\CC_{n,k}}\E_{(m,n,p_{x,n},q_{x,n})}\psi-\E_{(m,n,h_n,h_n)}\psi\nonumber\\
&\leq\ \frac{1}{\sharp\CC_{n,k}}\sum_{x\in\CC_{n,k}}\E_{(m,n,p_{x,n},q_{x,n})}\psi-\E_{(m,n,h_n,h_n)}\psi\nonumber\\ 
&\leq\ \E_{(m,n,h_n,h_n)}\bigg\arrowvert\frac{1}{\sharp \CC_{n,k}}\sum_{x\in\CC_{n,k}}\frac{d\P_{(m,n,p_{x,n},q_{x,n})}}{d\P_{(m,n,h_n,h_n)}}(\underbar{X})-1\bigg\arrowvert.\label{eq: likelihood}
\end{align}
For short we write $\E_0$ for $\E_{(m,n,h_n,h_n)}$ in the sequel. Note that the test is allowed to depend on the nuisance functional $h_n$ (in fact the $\log$-likelihood and its distribution do). Now we aim at determining $\delta_n$ such that the right-hand-side tends to zero as
$n$ goes to infinity. 
 Although $\lambda\big(\supp(\phi_{x,n})\cap\supp(\phi_{y,n})\big)=0$ for any different $x,y\in\CC_{n,k}$, the likelihood-ratios
\begin{align*}
L_{x,n}\ :=&\ \frac{d\P_{(m,n,p_{x,n},q_{x,n})}}{d\P_{(m,n,h_n,h_n)}}(\underbar{X})\ =\
\prod_{i=1}^{m}\biggl(1+\big(1-(m/n)\big)
\frac{\phi_{x,n}}{\sqrt{h_n}}(X_i)\biggr)\prod_{i=m+1}^n\biggl(1-
(m/n)\frac{\phi_{x,n}}{\sqrt{h_n}}(X_i)\biggr),
\end{align*}
are not independent. However, they are independent conditional on the random vector $\Delta_n=(\Delta_{x,n})_{x\in\CC_{k,n}}$ with entries
$$
  \Delta_{x,n}\ :=\ \Big(\sharp\big\{i\leq m: \Arrowvert X_i-x\Arrowvert_2\leq k\delta_n\ \big\},\ \sharp\big\{i> m: \Arrowvert X_i-x\Arrowvert_2\leq k\delta_n\ \big\}\Big).
$$ 
Note that $\E_0(L_{x,n}\arrowvert\Delta_n)=\E_0\, L_{x,n}=1$. 
Following at this point standard truncation arguments as, for instance, in D\"umbgen and Walther (2009), proof of Lemma 7.4, it turns out to be sufficient
for the convergence to zero of (\ref{eq: likelihood}) to find $\delta_n$ and $\gamma=\gamma_n\in (0,1]$ such that the ratio
\begin{equation}\label{eq: ee}
\max_{x\in\CC_{n,k}}\frac{1}{(\sharp\,\CC_{n,k})^{\gamma}}\E_0\,
L_{x,n}^{1+\gamma}
\end{equation}
tends to zero as $n$ goes to infinity. But
\begin{align}
 \E_0 L_{x,n}^{1+\gamma} &= \bigg\{\int
h_n(z)\biggl(1+(1-m/n)\frac{\phi_{x,n}(z)}{\sqrt{h_n(z)}}\biggr)^{1+\gamma}dz\biggr\}^{m}\bigg\{\int
h_n(z)\biggl(1-(m/n)\frac{\phi_{x,n}(z)}{\sqrt{h_n(z)}}\biggr)^{1+\gamma}dz\biggr\}^{n-m}\nonumber\\
&=
\bigg\{1+\frac{1}{2}\gamma(1+\gamma)\Big(1+O\big(\delta_n^{\beta}\big)\Big)(1-(m/n))^2\int_0^1\phi_{x,n}(z)^2dz\bigg\}^{m}\times
\label{eq: 1}\\
&\ \ \ \ \ \ \ \ \ \ \ \ \ \ \ \ \ \ \ \ \ \ \ \ \ \ \ \ \ \ \ 
\bigg\{1+\frac{1}{2}\gamma(1+\gamma)\Big(1+O\big(\delta_n^{\beta}\big)\Big)(m/n)^2\int_0^1\phi_{x,n}(z)^2dz\bigg\}^{n-m},\nonumber
\end{align}
using the bound 
$
(1+\Delta)^{1+\gamma} \leq 1 + (1+\gamma)\Delta +
2^{-1}\gamma(1+\gamma)\Delta^2 +3\gamma\Delta^2\arrowvert
\Delta\arrowvert
$
for $\arrowvert \Delta\arrowvert\leq 1$. Now let $\tilde{\phi}_k$ be the solution to the following optimization problem

\vspace{1mm}
\noindent
($**$) Minimize $\Arrowvert g\Arrowvert_2$ subject to
\begin{equation}\label{eq: ORR}
g\in\HH_d(\beta,L;\R^d),\ \ \supp(g)\subseteq B_0(k),\ \ g(0)=1,\ \ \int g(x)dx =0.
\end{equation}
Notice the rescaling property $L\delta_n^{\beta}g(./\delta_n)\in\HH_d(\beta,L;\R^d)$ with $\supp\big(L\delta_n^{\beta}g(./\delta_n)\big)=B_0(\delta_nk)$ and $L\delta_n^{\beta}g(0)=L\delta_n^{\beta}$ $\Leftrightarrow$ $g\in\HH_d(\beta,L;\R^d)$ with $\supp(g)=B_0(k)$ and $g(0)=1$.
Due to the equicontinuity of $(h_n)_{n\in\N}$, 
$$
\lim_{\delta\searrow 0}\sup_{x\in B_{z}(\delta)}\sup_n\big\arrowvert h_n(x)-h_n(z)\big\arrowvert =0,
$$
whence
\begin{equation}
\int \phi_{x,n}(z)^2\,dz\  
=\ 
\bigl(1+o(1)\bigr)L^2\delta_n^{2\beta+d}\Arrowvert\tilde{\phi}_k\Arrowvert_2^2
\label{eq: 2}
\end{equation}
because the minimum in ($*$) depends continuously on the mixed density $h_n$ as can be seen using a Lagrange multiplier for the centering constraint. Note that the $o(1)$-term is uniformly in $x\in\CC_{k,n}$. 
Now the combination of (\ref{eq: 1}) and (\ref{eq: 2}) shows that for $\delta_n$
sufficiently small, (\ref{eq: ee}) is bounded by
$$
\exp\biggl(n(m/n)(1-m/n)\frac{1}{2}\gamma(1+\gamma)L^2\delta_n^{2\beta+d}\Arrowvert
\tilde{\phi}_k\Arrowvert_2^2\big(1+o(1)\big)-\gamma\log(\sharp\, \CC_{k,n})\biggr).
$$
By construction, $\sharp\,\CC_{k,n} \geq d_k\cdot\delta_n^{-d}$ for some constant $d_k>0$. Now fix $\delta>0$ and define
$$
c_k(\beta,L)\ :=\ \biggl(\frac{2\,d\,L^{d/\beta}}{(2\beta
  +d)\Arrowvert
  \tilde{\phi}_k\Arrowvert_2^2}\biggr)^{\beta/(2\beta+d)}.
$$
Observe that the sequence $c_k(\beta,L)$ is increasing in $k$. We need to check  that 
$
\lim_{k\rightarrow\infty}\Arrowvert \tilde{\phi}_k\Arrowvert_2 =\Arrowvert \gamma_{\beta}\Arrowvert_2
$. Note that in contrast to (\ref{eq: ORR}), the solution of (\ref{eq: recovery}) does not integrate to zero in general and it remains still open if $\gamma_{\beta}$ is compactly supported for $d\geq 2$ and $\beta>1$. Starting from $\gamma_{\beta}$, it is sufficient to construct a sequence $\tilde{\gamma}_{\beta,k}$ satisfying the constraints of the optimization problem ($**$) such that $\lim_{k\rightarrow\infty}\Arrowvert \tilde{\gamma}_{\beta,k}\Arrowvert_2=\Arrowvert \gamma_{\beta}\Arrowvert_2$. Then the equality $
\lim_{k\rightarrow\infty}\Arrowvert \tilde{\phi}_k\Arrowvert_2 =\Arrowvert \gamma_{\beta}\Arrowvert_2$  follows from $\Arrowvert\tilde{\gamma}_{\beta,k}\Arrowvert_2\geq \Arrowvert\tilde{\phi}_k\Arrowvert_2$. The existence is sketched in the appendix of the extended version of this article. 
As a consequence there exists some $k'\in \N$ such that $c(\beta,L)(1-\delta) < c_{k'}(\beta,L)(1-\delta/2)$. Now one verifies that the lower bound is established with the choice 
$$
\delta_n\ :=\ \Bigl(\frac{c_{k'}(\beta,L)(1-\delta/2)\rho_n}{L}\Big)^{1/\beta}
$$  
and some sequence $\gamma=\gamma_n\rightarrow 0$ with $\lim_n\gamma_n(\log n)^{1/2}=\infty$.
\hfill$\square$

\medskip
\paragraph{\sc Proof of Theorem \ref{thm: Effizienz2}}  
By virtue of Theorem \ref{Thm: limit distribution}, the sequence $\LL\big(T_n\circ\Pi\big\arrowvert \XX_n\big)$ is tight in $\big(\P_n^{\otimes m}\otimes \Q_n^{\otimes (n-m)}\big)$-probability, resulting in stochastic boundedness of the sequence of random quantiles $\big(\kappa_{\alpha}(\underbar{X})\big)_{n\in\N}$. The bounded total variation of the kernel for $\beta\leq1$ is a consequence of its monotonicity, for $\beta > 1$ it results from the continuous differentiability of $\psi_{\beta,K}$ and its compact support. 
For notational convenience the dependency on $\beta$ and $K$ is suppressed. They are arbitrary but fixed unless stated otherwise. First note that for any  random couple $(\hat{j}_n,\hat{k}_n)$ it holds true that
\begin{align*}
\P_{(m,n,p_n,q_n)}\Big(T_n > \kappa_{\alpha}(\underbar{X})\Big)\ \geq\ \P_{(m,n,p_n,q_n)}\Big(T_{\hat{j}_n\hat{k}_nn}-C_{\hat{j}_n\hat{k}_nn}\, > \, \kappa_{\alpha}(\underbar{X})\Big).
\end{align*}
Hence it is sufficient to prove that for any sequence $(\phi_n)_{n\in\N}$ of admissible alternatives there exists a random sequence of $(\hat{j}_n,\hat{k}_n)_{n\in\N}$ with 
$
T_{\hat{j}_n\hat{k}_nn} - C_{\hat{j}_n\hat{k}_nn} \longrightarrow_{\P^{\otimes m}\otimes\Q^{\otimes (n-m)}}\infty.
$
As in the proof of Theorem \ref{Thm: limit distribution} define $\gamma_n(t,r)^2:=\E\hat{\gamma}_{2,n}(t,r)^2-\big(\E\hat{\gamma}_{1,n}(t,r)\big)^2$, $(t,r)\in\TT$.  Let 
$
t_n:=\argmax_{x\in J}\arrowvert \phi_n(x)\arrowvert$ and $r_n:=\big(\Arrowvert\phi_n\Arrowvert_{\sup}\big/L\big)^{1/\beta}$.
Define
$
\big(\hat{t}_n,\hat{r}_n\big):=\big(X_{\hat{j}_n},\big\Arrowvert X_{\hat{j}_n}-X_{\hat{k}_n}\big\Arrowvert_2\big)$ with
$$
(\hat{j}_n,\hat{k}_n)\ :=\ \argmin_{j,k=1,...,n}\lambda\Big(B_{t_n}(r_n)\,\Delta\, B_{X_j}\big(\Arrowvert X_j-X_k\Arrowvert_2\big)\Big).
$$
Now let the process $S_n$ on $\TT$ pointwise be defined by
$$
S_n\big(t,r\big)\ :=\ \frac{\sqrt{\lambda_n(1-\lambda_n)}}{\sqrt{n}}\sum_{i=1}^n\psi\Big(\frac{\Arrowvert X_i-t\Arrowvert_2}{r}\Big)\Lambda(X_i).
$$

\vspace{-2mm}
\noindent
Furthermore, let us introduce the random variables $(\hat{t}_{ni},\hat{r}_{ni})$, based on the indices $(\hat{j}_{ni},\hat{k}_{ni})$ which are defined analogously to $(\hat{j}_n,\hat{k}_n)$ but with the minimum running over the set $j,k\in\{1,...,n\}\setminus\{i\}$ only.
Then, recalling the definition $
\psi_{tr}(x)\ :=\ \psi\big(\frac{\Arrowvert t-x\Arrowvert_2}{r}\big)
$,
\begin{align}
\frac{1}{\gamma_n(t_n,r_n)}&\Big\arrowvert \E\Big(S_n(\hat{t}_n,\hat{r}_n)-S_n(t_n,r_n)\Big)\Big\arrowvert\nonumber\\
&=\ \frac{\sqrt{\lambda_n(1-\lambda_n)}}{\gamma_n(t_n,r_n)}\frac{1}{\sqrt{n}}\bigg\arrowvert \frac{n}{m}\sum_{i=1}^m\E\Big(\psi_{\hat{t}_n\hat{r}_n}(X_i)-\psi_{t_nr_n}(X_i)\Big)\nonumber\\ &\ \ \ \ \ \ \ \ \ \ \ \ \ \ \ \ \ \ \ \ \ \ \ \ \ \ -\  \frac{n}{n-m}\sum_{i=m+1}^n\E\Big(\psi_{\hat{t}_n\hat{r}_n}(X_i)-\psi_{t_nr_n}(X_i)\Big)\bigg\arrowvert\nonumber\\
&\leq\  \frac{\sqrt{\lambda_n(1-\lambda_n)}}{\gamma_n(t_n,r_n)}\frac{1}{\sqrt{n}}\bigg\arrowvert \frac{n}{m}\sum_{i=1}^m\E\Big(\psi_{\hat{t}_{n}\hat{r}_{n}}(X_i)-\psi_{\hat{t}_{ni}\hat{r}_{ni}}(X_i)\Big)\nonumber\\ 
&\ \ \ \ \ \ \ \ \ \ \ \ \ \ \ \ \ \ \ \ \ \ \ \ \ \ -\ \frac{n}{n-m}\sum_{i=m+1}^n\E\Big(\psi_{\hat{t}_{n}\hat{r}_{n}}(X_i)-\psi_{\hat{t}_{ni}\hat{r}_{ni}}(X_i)\Big)\bigg\arrowvert\nonumber\\
&\ \ \ \ \ +\ \frac{\sqrt{\lambda_n(1-\lambda_n)}}{\gamma_n(t_n,r_n)}\frac{1}{\sqrt{n}}\bigg\arrowvert \frac{n}{m}\sum_{i=1}^m\E\Big(\psi_{\hat{t}_{ni}\hat{r}_{ni}}(X_i)-\psi_{t_nr_n}(X_i)\Big)\nonumber\\ 
&\ \ \ \ \ \ \ \ \ \ \ \ \ \ \ \ \ \ \ \ \ \ \ \ \ \ \ \ \ \ \ -\ \frac{n}{n-m}\sum_{i=m+1}^n\E\Big(\psi_{\hat{t}_{ni}\hat{r}_{ni}}(X_i)\, -\, \psi_{t_nr_n}(X_i)\Big)\bigg\arrowvert\nonumber\\
&\leq\   \frac{\sqrt{\lambda_n(1-\lambda_n)}}{\gamma_n(t_nr_n)}\frac{4}{\sqrt{n}}\Arrowvert \psi\Arrowvert_{\sup}\max\Big(\frac{n}{m},\frac{n}{n-m}\Big)\nonumber\\
&\ \ \ \ \  +\
 \frac{\sqrt{\lambda_n(1-\lambda_n)}}{\gamma_n(t_n,r_n)}\frac{1}{\sqrt{n}}\Bigg\arrowvert\E\,\Bigg\{ \frac{n}{m}\sum_{i=1}^m\int\Big(\psi_{\hat{t}_{ni}\hat{r}_{ni}}(x)-\psi_{t_nr_n}(x)\Big)p_n(x)dx\label{expr}\\
&\ \ \ \ \ \ \ \ \ \ \ \ \ \ \ \ \ \ \ \ \ \ \ \ \ \ \ \ \ \ \ \ \ -\ \frac{n}{n-m}\sum_{i=m+1}^n\int\Big(\psi_{\hat{t}_{ni}\hat{r}_{ni}}(x)-\psi_{t_nr_n}(x)\Big)q_n(x)dx\Bigg\}\Bigg\arrowvert,\nonumber
\end{align}
whereby we used for the first term in the last inequality that $(\hat{t}_{ni},\hat{r}_{ni})$ differs from $(\hat{t}_n,\hat{r}_n)$ for at most two indices $i,j\in\{1,...,n\}$; the second term follows by including and evaluating the conditional expectation given $(\hat{t}_{ni},\hat{r}_{ni})$ as $X_i$ is independent of $(\hat{t}_{ni},\hat{r}_{ni})$. Replacing again $(\hat{t}_{ni},\hat{r}_{ni})$ by $(\hat{t}_n,\hat{r}_n)$, the second expression behind the inequality in formula (\ref{expr}) is bounded by
\begin{align}
 \frac{\sqrt{\lambda_n(1-\lambda_n)}}{\gamma_n(t_n,r_n)}\frac{4}{\sqrt{n}}&\Arrowvert \psi\Arrowvert_{\sup}\max\Big(\frac{n}{m},\frac{n}{n-m}\Big)\nonumber\\ 
& +\ \frac{\sqrt{n}\sqrt{\lambda_n(1-\lambda_n)}}{\gamma_n(t_n,r_n)}\Bigg\arrowvert\E\bigg[\int\Big(\psi_{\hat{t}_n\hat{r}_n}(x)-\psi_{t_nr_n}(x)\Big)\Big(p_n(x)-q_n(x)\Big)dx\bigg]\Bigg\arrowvert. \label{expr2}
\end{align}
Now we can make use of the fact that
$
\big\arrowvert p_n(x)-q_n(x)\big\arrowvert = \big\arrowvert \phi_n(x)\sqrt{h_n(x)}\big\arrowvert \leq C\Arrowvert \phi_n\Arrowvert_{\sup}
$
with $C:=\sup_n\sup_x\big\arrowvert \sqrt{h_n(x)}\big\arrowvert$. Recall that $\Arrowvert h_n\Arrowvert_{\sup}$ is uniformly bounded due to the equicontinuity assumption on $(h_n)_{n\in\N}$ and the constraint on the $L_1$-nor\def\N{\mathbb{N}}m $\Arrowvert h_n\Arrowvert_1=1$, 
whence 
the term in (\ref{expr2}) is not greater than
\begin{equation}
C\frac{\sqrt{n}\Arrowvert \phi_n\Arrowvert_{\sup}}{\gamma_n(t_n,r_n)}\,\E\bigg(\int\Big\arrowvert\psi_{\hat{t}_n\hat{r}_n}(x)-\psi_{t_nr_n}(x)\Big\arrowvert dx\bigg).\label{expr3}
\end{equation}
Using the bounded total variation $TV(\psi)$ of $\psi$ and $
M_x
$ and $\mu$ as defined in the proof of Theorem \ref{Thm: limit distribution},  
the integral which appears in (\ref{expr3}) can be bounded by
\begin{align}
  \E\bigg(\int\Big\arrowvert\psi_{\hat{t}_{n}\hat{r}_{n}}(x)&-\psi_{t_nr_n}(x)\Big\arrowvert dx\bigg)\nonumber\\
&\leq\   \E\bigg(\int \mu\big( M_x(t_n,r_n,\hat{t}_{n},\hat{r}_{n})\big)dx\bigg)\nonumber\\
&=\ \E\bigg(\int\int I\big\{y\in M_x(t_n,r_n,\hat{t}_{n},\hat{r}_{n})\big\}dxd\mu(y)\bigg)\ \ \ \ \ \ \ \ \  \text{(Fubini)}\nonumber\\
&\leq\   TV(\psi)\,\E\sup_{y\in [0,1]}\bigg(\int  I\big\{y\in M_x(t_n,r_n,\hat{t}_{n},\hat{r}_{n})\big\}dx\bigg)\nonumber\\
&= \  TV(\psi)\,\E\,\lambda\Big(B_{t_n}(r_n)\Delta B_{\hat{t}_{n}}(\hat{r}_{n})\Big)\nonumber\\
&=\ O\big(r_n^{d-1}n^{-1/d}\big),\label{eq: rate}
\end{align}
using in the last bound besides the stochastic convergence rate $n^{-1/d}$ the uniform integrability of the sequences $\big(n^{1/d}\Arrowvert \hat{t}_n-t_n\Arrowvert_2\big)$, $\big(n^{1/d}\arrowvert\hat{r}_n-r_n\arrowvert\big)$ which result from $\P\big(\Arrowvert \hat{t}_{n}-t_n\Arrowvert_2>x\big)=\prod_{i=1}^n\P\big(X_i\not\in B_{t_n}(x)\big)\sim \big(1-\lambda(B_{t_n}(x)\cap[0,1]^d)\big)^{n}\ (=(1-Vx^d)^n\text{~if~} B_{t_n}(x)\subset[0,1]^d)$ and $\P\big(\arrowvert \hat{r}_{n}-r_n\arrowvert>x\big)\leq 2\,\P\big(\Arrowvert \hat{t}_{n}-t_n\Arrowvert_2>x/2\big)$. Here, $V$ denotes the volume of the $d$-dimensional Euclidean unit ball, i.e. $V=\pi^{d/2}\Gamma(d/2+1)$. Together with (\ref{expr}) - (\ref{expr3}) this shows that for any sequence of admissible alternatives $(\phi_n)_{n\in\N}$
\begin{equation}
\frac{\big\arrowvert\E\big( {S}_n(\hat{t}_n,\hat{r}_n)-{S}_n(t_n,r_n)\big)\big\arrowvert}{\gamma_n(t_n,r_n)}\ =\ O\Big(r_n^{d/2-1+\beta}n^{-1/d+1/2}\Big).\label{eq: e-approximation}
\end{equation}
If in particular $\Arrowvert \phi_n\Arrowvert_{\sup}=O\Big(\big((\log n)/n\big)^{\beta/(2\beta+d)}\Big)$, the term in (\ref{eq: e-approximation}) is of order
$$
O\Big((\log n)^{(\beta+d/2-1)/(2\beta+d)}n^{-(2\beta/d)/(2\beta+d)}\Big).
$$
We need to check that
\begin{equation}
\frac{\gamma_n(t_n,r_n)}{\hat{\gamma}_n(\hat{t}_n,\hat{r}_n)}\, \longrightarrow_{\P^{\otimes m}\otimes Q^{\otimes (n-m)}}\ 1.\label{eq: convergence variance}
\end{equation}
For this we use the decomposition 
$
[(n-1)/n]\hat{\gamma}_n(t,r)^2\ =\ \hat{\gamma}_{2,n}(t,r)^2\ -\ \hat{\gamma}_{1,n}(t,r)^2
$
To this end note first that
\begin{align*}
\big\arrowvert \hat{\gamma}_{n,1}&(\hat{t}_n,\hat{r}_n)-\hat{\gamma}_{n,1}(t_n,r_n)\big\arrowvert \\ &\leq\ \big\Arrowvert \psi_{\hat{t}_n\hat{r}_n}-\psi_{t_nr_n}\big\Arrowvert_{\sup}\frac{1}{n}\sum_{i=1}^nI\Big\{X_i\in B_{\hat{t}_n}\big(\hat{r}_n\big)\cap B_{t_n}\big(r_n\big)\Big\}\\ 
&\ \ \ \ \ +\ 2\Arrowvert \psi\Arrowvert_{\sup}\frac{1}{n}\sum_{i=1}^nI\Big\{X_i\in B_{\hat{t}_n}\big(\hat{r}_n\big)\Delta B_{t_n}\big(r_n\big)\Big\}\\
&\leq\ \big\Arrowvert \psi_{\hat{t}_n\hat{r}_n}-\psi_{t_nr_n}\big\Arrowvert_{\sup}\frac{1}{n}\sum_{i=1}^nI\Big\{X_i\in B_{t_n}\big(r_n\big)\Big\}\\ 
&\ \ \ \ \ +\ 2\Arrowvert \psi\Arrowvert_{\sup}\frac{1}{n}\sum_{i=1}^nI\Big\{X_i\in B_{\hat{t}_n}\big(\hat{r}_n\big)\Delta B_{t_n}\big(r_n\big)\Big\}\\
&=\ o_p(1)O_p(r_n^d)\ +\ O_p\big(r_n^{d-1}n^{-1/d}\big)\ =\ o_p\Big(\gamma_{n,1}(t_n,r_n)\Big).
\end{align*} 
The "$o_p(1)$"-term results from the H\"older continuity of $\psi$ (for $\beta>1$ the first derivative of $\psi$ is uniformly bounded on $[-K,K]$), $\supp\big(\psi_{t_nr_n}-\psi_{\hat{t}_n\hat{r}_n}\big)=B_{t_n}(r_n)\cup B_{\hat{t}_n}(\hat{r}_n)$ and the fact that $r_n >\big(c(\beta,L)\rho_{m,n}/L\big)^{1/\beta}$ while $\hat{t}_n-t_n\sim n^{-1/d},\ \hat{r}_n-r_n\sim n^{-1/d}$. 
The case $i=2$ is done analogously (taking the square). To verify (\ref{eq: convergence variance}) it remains to be shown that
$
\hat{\gamma}_n({t}_n,{r}_n)/\gamma_n(t_n,r_n) -1\, =\, o_p( 1)
$
which however is a simple consequence of Chebychef's inequality since for any $\beta>0$ and any sequence of admissible alternatives $(\phi_n)_{n\in\N}$, the sequence $\gamma_n(t_n,r_n)\sim r_n^{d/2}$ or some subsequence decreases (if it decreases) at a slower rate than $n^{-1/2}$.
The above considerations show in particular that 
\begin{align*}
C_{\hat{j}_n\hat{k}_nn}\ &=\ \frac{3\,R_{\psi}(m,n)}{\sqrt{n}\,\hat{\gamma}_n(\hat{t}_n,\hat{r}_n)}\delta(m,n)\log\Big(\hat{\gamma}_n(\hat{t}_n,\hat{r}_n)^{-2}\Big)\ +\ \delta(m,n)\sqrt{2\log\Big(\hat{\gamma}_n(\hat{t}_n,\hat{r}_n)^{-2}\Big)}\\
&=\ \sqrt{2\log\Big({\gamma}_n({t}_n,{r}_n)^{-2}\Big)}\ +\ o_p(1),
\end{align*}
using in addition that $\delta(m,n)=1+O(n^{-1/2})$. Consequently,
\begin{align}\label{term}
T_{\hat{j}_n\hat{k}_nn}\,-\, C_{\hat{j}_n\hat{k}_nn}\ &=\ O_p(1)\  +\ \frac{\E {S}_n(t_n,r_n)}{\gamma_n(t_n,r_n)}\Big(1\, +\, o_p(1)\Big)\ -\ \sqrt{2\log\Big({\gamma}_n({t}_n,{r}_n)^{-2}\Big)},
\end{align}
and it has to be verified that the latter quantity goes to infinity. Recall that
\begin{align}
\gamma_n(t_n,r_n)^2\ &=\ \int_{[0,1]^d} \psi_{t_nr_n}(x)^2h_n(x)dx\ -\ \bigg(\int_{[0,1]^d}\psi_{t_nr_n}(x)h_n(x)dx\bigg)^2\nonumber\\
&=\ \Big(1+O(r_n^d)\Big)\int_{[0,1]^d} \psi_{t_nr_n}(x)^2h_n(x)dx.\label{appp}
\end{align}
We first assume that $r_n=o(1)$, i.e. $\Arrowvert\phi_n\Arrowvert_{\sup}=o(1)$. 
Using that
$$
\lim_{\delta\searrow 0}\sup_n\sup_{t\in[0,1]^d}\sup_{x\in B_{t}(\delta)}\big\arrowvert h_n(x)-h_n(t)\big\arrowvert\ =\ 0,
$$
which follows by  the same argument as used in Theorem \ref{Thm: lower bound} and the fact that any sequence of centers $(t_n)_{n\in\N}$ has a convergent subsequence by the compactness
of $[0,1]^d$, 
\begin{equation}\label{eq: 1pppp}
\frac{\E {S}_n(t_n,r_n)}{\gamma_n(t_n,r_n)}\ =\ \sqrt{n}\sqrt{\lambda_n(1-\lambda_n)}\frac{\int_{[0,1]^d}\psi_{t_nr_n}(x)\phi_n(x)dx}{\big[\int_{[0,1]^d}\psi_{t_nr_n}(x)^2dx\big]^{1/2}}\,\big(1+o(1)\big).
\end{equation}
Using the approximation in (\ref{appp}) we obtain analogously
\begin{equation}\label{eq: apppp}
\sqrt{2\log\Big({\gamma}_n({t}_n,{r}_n)^{-2}\Big)}\ =\ \bigg[2\log\bigg(1\bigg/O(1)\int_{[0,1]^d}\psi_{t_nr_n}(x)^2dx\bigg)\bigg]^{1/2}.
\end{equation}
Recall that $\psi=\psi_{\beta,K}$ with $K$ the bound of the support. Standard calculation shows that the bounded $L_2$-norm of $\gamma_{\beta}$ implies 
$$
\frac{\big\arrowvert\int\psi_{t_nr_n;\beta,K}(x)\phi_n(x)dx\,\big\arrowvert}{\big[\int\psi_{t_nr_n;\beta,K}(x)^2dx\big]^{1/2}}\ =\ \frac{\big\arrowvert\int\psi_{t_nr_n;\beta}(x)\phi_n(x)dx\,\big\arrowvert}{\big[\int\psi_{t_nr_n;\beta}(x)^2dx\big]^{1/2}}\big(1+c_K\big)\ \ \text{with $c_K\rightarrow 0$ as~$K\rightarrow\infty$,}
$$
but note that the total variation $TV(\psi_{\beta,K})$ is increasing in $K$.
Define now $\delta_n:=(1+\delta)c(\beta,L)\rho_{m,n}$. Then by its construction, $\delta_n\psi_{t_nr_n;\beta}\in\HH_d\big(\beta,L;\R^d\big)$. Moreover, by the closedness in $L_2$ and the convexity of the sets $\big\{\phi\in\HH_d(\beta,L;\R^d): \phi(t_n)\geq \delta_n\big\}$ and $\big\{\phi\in\HH_d(\beta,L;\R^d): \phi(t_n)\leq -\delta_n\big\}$, it results finally from convex analysis and the definition of $\gamma_{\beta}$ that 
$$
\frac{\big\arrowvert\int\psi_{t_nr_n;\beta}(x)\phi_n(x)dx\,\big\arrowvert}{\big[\int\psi_{t_nr_n;\beta}(x)^2dx\big]^{1/2}}\ \geq\ \frac{\delta_n^{-1}\Arrowvert \delta_n\psi_{t_nr_n;\beta}\Arrowvert_2^2}{\Arrowvert\psi_{t_nr_n;\beta}\Arrowvert_2}\ =\ \delta_nr_n^{d/2}\Arrowvert\gamma_{\beta}\Arrowvert_2.
$$
Combining (\ref{appp}) -- (\ref{eq: apppp}), one verifies for the expression of the right hand side in (\ref{term}) that it possesses the approximation
\begin{align*}
(\ref{term})\ &=\ O_p(1)\ +\ \sqrt{n}\sqrt{\lambda_n(1-\lambda_n)}\delta_n r_n^{d/2}\Arrowvert\gamma_{\beta}\Arrowvert_2\big(1+c_K\big)\ -\ \Big(\frac{2d}{2\beta+d}\Big)^{1/2}\sqrt{\log\big(n\big/\log n\big)}\\
&=\ O_p(1) \ +\ \sqrt{\log n}\bigg(\frac{2dL^{d/\beta}}{(2\beta+d)\Arrowvert \gamma_{\beta}\Arrowvert_2^2}\bigg)^{1/2}L^{-d/(2\beta)}\Arrowvert\gamma_{\beta}\Arrowvert_2(1+c_K)(1+\delta)^{d/(2\beta)+1}\\
&\ \ \ \ \ \ \ \ \ \ \ \ \ \ \ \ \ \ \  - \ \Big(\frac{2d}{2\beta+d}\Big)^{1/2}\sqrt{\log\big(n\big/\log n\big)},
\end{align*}
which goes to infinity for $K$ sufficiently large. If there exists a sequence $(\phi_n)_{n\in\N}$ of admissible alternatives such that
$
\limsup_{n\rightarrow\infty}\P_{(m,n,p_n,q_n)}\big(T_n >  \kappa_{\alpha}(\underbar{X})\big)\, <\, 1,
$
there exists by the considerations above a subsequence (for simplicity also denoted by $(n)$) along which $\Arrowvert \phi_{n}\Arrowvert_{\sup}$ stays uniformly bounded away from zero.  But the bounds (\ref{eq: e-approximation}) and (\ref{eq: convergence variance}) show that
$$
\frac{\E S_n(\hat{t}_n,\hat{r}_n)-\E S_n(t_n,r_n)}{\gamma_n(\hat{t}_n,\hat{r}_n)}\ =\ O\Big(n^{-1/d+1/2}\Big)\big(1+o_p(1)\big),
$$
as well as the logarithmic correction term $C_{\hat{j}_n\hat{k}_nn}$ are in this case of smaller order than $\arrowvert\E S_n(t_n,r_n)\arrowvert$, which concludes the proof by contradiction. \hfill$\square$

\medskip
\paragraph{\sc Proof of Theorem \ref{thm: spatial adaptivity}}
Following the considerations of the proof of Theorem \ref{thm: Effizienz2}, it has to be established that there exist random sequences $\big(\hat{j}_{ni},\hat{k}_{ni}\big)_{n\in\N}$ with $B_{X_{\footnotesize\hat{j}_{ni}}}\big(\big\Arrowvert X_{\hat{j}_{ni}}-X_{\hat{k}_{ni}}\big\Arrowvert_2\big)\,\subset\, J_i$, $i=1,...,k$, such that for any sequence of alternatives as formulated in Theorem \ref{thm: spatial adaptivity} and any fixed $K>0$
$$
\liminf_{n\rightarrow\infty}\ \P_{(m,n,p_n,q_n)}\Big(T_{\hat{j}_{ni}\hat{k}_{ni}}-C_{\hat{j}_{ni}\hat{k}_{ni}}>\kappa_{\alpha}(\underbar{X})\Big)\ =\ 1,\ \ i=1,...,k.
$$
Then the result follows because the finite intersection of sets with asymptotic probability equal to $1$ has asymptotically mass $1$ as well. Inspired by the arguments in Rohde (2008) for the univariate regression context, we first establish the following:

\medskip
For $\phi_n\in\HH_d\big(\beta,L;[0,1]^d\big)$ with $\Arrowvert \phi_n\Arrowvert_{\sup}\leq 1$ and $x^*=\argmax_{x\in [0,1]^d}\arrowvert \phi_n(x)\arrowvert$, there exists some constant $c=c(\beta, L)>0$ and a compact ball $B=B(\phi_n)\subset \R^d$ with center $x^*$ such that 
\begin{equation}\label{eq: function}
\lambda\big( B\cap [0,1]^d\big) \geq\  c\arrowvert \phi_n(x^*)\arrowvert^{d/\beta} \text{~~and~~} \big\arrowvert\phi_n(x)\big\arrowvert\ \geq\ \frac{1}{2}\big\arrowvert \phi_n(x^*)\big\arrowvert\ \ \text{for all $x\in B\cap[0,1]^d$}.
\end{equation}
Assume that $\beta>1$ (the above inequality is trivial in case $\beta\leq 1$). With $j=(j_1,...,j_d)$ we denote subsequently some multi-index, where $\arrowvert j\arrowvert = j_1+...+j_d$ defines its length, $x^j:=\prod_{i=1}^dx_i^{j_i}$ and $D^j:=\partial^{\arrowvert j\arrowvert}\big/[\partial x_1^{j_1}\cdot ...\cdot \partial x_m^{j_m}]$ the partial differential operator. 
Let $\phi\in\HH_d\big(\beta,L;[0,1]^d\big)$ with $\Arrowvert\phi\Arrowvert_{\sup}=D>0$. By the definition of the isotropic H\"older class we have $\arrowvert\phi(x)-T_y^{(f)}(x)\arrowvert\leq L\Arrowvert x-y\Arrowvert_2^{\beta}\ (\leq L\sqrt{d}^{\beta})$, which entails that $\sup_y \Arrowvert T_y^{(f)}\Arrowvert_{[0,1]^d}\leq D+L\sqrt{d}^{\beta}$. In order to establish (\ref{eq: function}), note that for any polynomial $P=\sum_{\arrowvert j\arrowvert\leq\lfloor \beta\rfloor}a_jx^j$, the topology induced by the metrics corresponding to the two norms $\Arrowvert P\Arrowvert_{(1)}=\sup_{x\in[0,1]^d}\arrowvert P(x)\arrowvert$ and $\Arrowvert P\Arrowvert_{(2)}:=\max_j\arrowvert a_j\arrowvert$ respectively on the ring of polynomials of total degree at most $\lfloor \beta\rfloor$ on $[0,1]^d$ is the topology of uniform convergence, hence these two norms are equivalent. Consequently, the boundedness of the polynomial $T_y^{(f)}$ by $D+L{\sqrt{d}~}^{\beta}$ uniformly in $y$ implies that there exists some constant $C=C(\beta)$ such that $\Arrowvert D^j\phi\Arrowvert_{\sup}\leq C\big(D+L)$ for all multi-indices $j$ with $\arrowvert j\arrowvert\leq\lfloor\beta\rfloor$. Now the Mean Value Theorem implies for some intermediate point $z\in \big\{x+t(x^*-x); 0\leq t\leq 1\big\}$ 
\begin{align*}
\big\arrowvert\phi(x)-\phi(x^*)\big\arrowvert \ &=\ \big\arrowvert\big( \nabla\phi(z)\big)^T\big(x-x^*\big)\big\arrowvert\\ 
&\leq\ \sqrt{d}\sup_{j:\, \arrowvert j\arrowvert=1}\big\Arrowvert D^j\phi\big\Arrowvert_{\sup}\Arrowvert x-x^*\Arrowvert_2\\ 
&\leq\ \sqrt{d}\,C\,\big(D+L\big)\Arrowvert x-x^*\Arrowvert_2.
\end{align*}
Thus, 
$$
\arrowvert \phi(x)\arrowvert\ \geq\ \frac{1}{2}\arrowvert \phi(x^*)\arrowvert\ \ \text{for all $x$ in}\ B_{x^*}\bigg(\frac{D}{2\sqrt{d}C(D+L)}\bigg)\cap[0,1]^d. 
$$
If $\phi\in\HH_d\big(\beta,L;[0,1]^d\big)$ with $\Arrowvert \phi\Arrowvert_{\sup}=\delta\leq 1$, then the function $g_{\delta}$, for $x\in[0,1]^d$ pointwise defined by $g_{\delta}(x):=\delta^{-1}\phi\big(\delta^	{1/\beta}x+x^*\big)\cdot I\big\{ \delta^{1/\beta}x+x^*\in[0,1]^d\big\}$ is element of $\HH_d(\beta,L; \supp (g_{\delta}))$ with $\Arrowvert g_{\delta}\Arrowvert_{\sup}=1$. Note that $\supp(g_{\delta})$ is a convex set. Therefore, the above considerations imply that $\arrowvert \phi(x)\arrowvert\geq \delta/2$ on
$$
B_{x^*}\bigg(\frac{\delta^{1/\beta}}{2\sqrt{d}C(1+L)}\bigg)\,\cap\, [0,1]^d.
$$
But then its Lebesgue measure is always greater than $c\arrowvert\delta\arrowvert^{d/\beta}$ for some constant $c=c(\beta,L)$, independent of $\delta$ and $x^*$, hence (\ref{eq: function}) is established.

\medskip
Let now $\beta_i,L_i\in (0,\infty)$ fixed but arbitrary, $J_i\subset [0,1]^d$ some nondegenerate rectangle, $\tilde{\phi}_n =p_n-q_n$ a sequence of functions with ${\tilde{\phi}}_{n\arrowvert J_i}\in\HH_d\big(\beta_i,L_i;J_i\big)$. It has to be shown that there exists a universal constant $k_i=k_i(\beta_i,L_i,c)$ such that $T_{\hat{j}_n\hat{k}_nn}-C_{\hat{j}_n\hat{k}_nn}\rightarrow_{\P^{\otimes m}\otimes\Q^{\otimes (n-m)}}\infty$ whenever 
$\Arrowvert \tilde{\phi}_n\Arrowvert_{J_i}\geq k_i\rho_{m,n}$. First, we choose a compact ball $B_i(\tilde{\phi}_n)$ with center $x_i^*:=\argmax_{t\in J_i}\arrowvert \tilde{\phi}_n(t)\arrowvert$ satisfying $\lambda(B_i(\tilde{\phi}_n)\cap J_i)\geq c\arrowvert \tilde{\phi}_n(x_i^*)\arrowvert^{d/\beta}$ and (\ref{eq: function}). Let the couple $(\hat{t}_n,\hat{r}_n):=\big(X_{\hat{j}_n},\Arrowvert X_{\hat{j}_n}-X_{\hat{k}_n}\Arrowvert_2\big)$ be defined by 
$$
(\hat{j}_n,\hat{k}_n)\ :=\ \argmin_{j,k\in\{1,...,n\}}\lambda\bigg(B_{X_j}\Big(\Arrowvert X_j-X_k\Arrowvert_2\Big)\,\Delta\, B_i(\tilde{\phi}_n)\bigg).
$$ 
Consulting the proof of Theorem \ref{thm: Effizienz2}, this definition of $(\hat{t}_n,\hat{r}_n)$ allows for an approximation as in (\ref{term}). Since $\arrowvert \tilde{\phi}_n(x)\arrowvert\geq 2^{-1}\Arrowvert \tilde{\phi}_n\Arrowvert_{J_i}$ for all $x\in B_i(\tilde{\phi}_n)\cap B_{\hat{t}_n}(\hat{r}_n)\cap J_i$,
\begin{align*}
\frac{\E S_n(t_n,r_n)}{\gamma_n(t_n,r_n)}\ &\geq\ \frac{1}{2}\Arrowvert\tilde{\phi}_n\Arrowvert_{J_i}\sqrt{n}\frac{\sqrt{\lambda_n(1-\lambda_n)}}{\sqrt{\max_x h_n(x)}}\Big[\E\lambda\Big(B_i(\tilde{\phi}_n)\cap B_{\hat{t}_n}(\hat{r}_n)\cap [0,1]^d\Big)\Big]^{1/2}\\ 
&\geq\  C\Arrowvert\tilde{\phi}_n\Arrowvert_{J_i}^{(\beta+d/2)/\beta}\sqrt{n}\big(1+o(1)\big) 
\end{align*}
for some universal constant $C=C\big(K,(\lambda_n)_{n\in\N}\big)>0$. Now the asserted result is easily deduced for a sufficiently large constant $k_i$.\hfill$\square$

\paragraph{\sc Proof of Theorem \ref{thm: local alternatives}.}

(i) Let $(p,q,m,n)$ be such that $h=h_n=I_{[0,1]^d}$ and $\phi_n$ the sequence of piecewise constant functions on $[0,1]^d$ with $\phi_n(z)=c_n/\sqrt{n\delta_n^d}$ for $z\in B_x(\delta_n)$, $\phi_n(z)=-c_n/\sqrt{n\delta_n^d}$ for $z\in B_x(\kappa\delta_n)\setminus B_{x}(\delta_n)$ and equals zero otherwise, where $\kappa=\kappa(d)>1$ is such that $\lambda\big(B_{x}(\kappa\delta_n)\setminus B_x(\delta_n)\big)=\lambda\big(B_x(\delta_n)\big)$
and $0<c_n\leq\sqrt{n\delta_n^d}$.
Then 
$$
\log\bigg(\frac{d\mathbb{P}_{(m,n,p_n,q_n)}}{d\mathbb{P}_{(m,n,h,h)}}(\underbar{X})\bigg)\ =\ \sum_{i=1}^m\log\Big(1+(1-m/n)\phi_n(X_i)\Big)\ +\ \sum_{j=m+1}^n\log\Big(1-(m/n)\phi_n(X_j)\Big)
$$
with $(X_k)_{k\in\N}$ iid uniformly distributed on $[0,1]^d$. Note that
\begin{align*}
\LL_{(m,n,h,h)}&\bigg[\log\bigg(\frac{d\mathbb{P}_{(m,n,p_n,q_n)}}{d\mathbb{P}_{(m,n,h,h)}}(\underbar{X})\bigg)\bigg]\\
 &\sim\ \sum_{i=1}^{N_m}\log\Big(1+(1-m/n)\frac{c_n}{\sqrt{n\delta_n^d}}\,R_i\Big)\ +\ \sum_{j=1}^{N_{n-m}}\log\Big(1-(m/n)\frac{c_n}{\sqrt{n\delta_n^d}}\,R_j'\Big)
\end{align*}
with $(R_{k})_{k\in\N}$ and $(R_{k}')_{k\in\N}$ two independent sequences of  iid Rademacher variables, $N_m$ and $N_{n-m}$ independent with
$$
N_m\sim\Bin\Big(m,V\kappa^d\delta_n^d\Big),\ \ N_{n-m}\sim\Bin\Big(n-m, V\kappa^d\delta_n^d\Big)\ \ \text{and}\ \ V\,=\,\pi^{d/2}\Gamma\big(d/2+1\big).
$$ 
Suppose first that $n\delta_n^d\not\rightarrow\infty$. By extracting a subsequence if necessary we may assume that $m/n\rightarrow\lambda\in(0,1)$, $c_n/\sqrt{n\delta_n^d}\rightarrow c\in[0,1]$ and $n\delta_n^d\rightarrow V^{-1}\kappa^{-d}\gamma$. Then, with $\rightarrow_w$ denoting weak convergence,  
\begin{equation}
 \LL_{(m,n,h,h)}\bigg[\log\bigg(\frac{d\mathbb{P}_{(m,n,p_n,q_n)}}{d\mathbb{P}_{(m,n,h,h)}}(\underbar{X})\bigg)\bigg] \longrightarrow_w\ \Q\label{eq: limit}
\end{equation}
with the convolution
$$
\mathbb{Q}\ :=\ \Bigg(\sum_{k=0}^{\infty}p_{\gamma(1-\lambda)}(k)\LL\bigg(\sum_{i=1}^{k}\log\Big(1-\lambda cR_i\Big)\bigg)\Bigg)\star\Bigg( \sum_{k'=0}^{\infty}p_{\gamma\lambda}(k')\LL\bigg(\sum_{j=1}^{k'}\log\Big(1+(1-\lambda)cR_j'\Big)\bigg)\Bigg)
$$
and the Poisson weights $p_{\mu}(k):=e^{-\mu}\mu^k/k!$. Since $\int e^zd\,\mathbb{Q}(z)=1$, we can apply Le Cam's notion of contiguity (Le Cam and Yang 2000, chapter 3) to conclude that
$$
\limsup_{n\rightarrow\infty}\, \mathbb{E}_{(m,n,p_n,q_n)}\psi_n(\underbar{X})\ <\ 1.
$$
Consequently $n\delta_n^d\rightarrow\infty$. Now assume that $n\delta_n^d\rightarrow\infty$ but $c_n\not\rightarrow\infty$. Without loss of generality we may assume that $c_n\rightarrow c'/\sqrt{V\kappa^d}\in[0,\infty)$. Then Lindeberg's  CLT entails that (\ref{eq: limit}) holds true with
$$
\mathbb{Q} \ :=\ \NN\bigg(-\frac{(1-\lambda)\lambda^2{c'}^2}{2},\, (1-\lambda)\lambda^2{c'}^2\bigg)\,\star\,\NN\bigg(-\frac{\lambda(1-\lambda)^2{c'}^2}{2},\, \lambda(1-\lambda)^2{c'}^2\bigg).
$$ 
Again, the limiting distribution satisfies $\int e^zd\,\mathbb{Q}(z)=1$, whence $c_n\rightarrow\infty$.

\medskip
(ii) We begin as in the proof of Theorem \ref{thm: Effizienz2}, but with $t_n:=x$, $r_n:=\delta_n$ and $\Arrowvert \phi_n\Arrowvert_{\sup}:=c_n/\sqrt{n\delta_n^d}$. Adjusting (\ref{expr}) -- (\ref{eq: e-approximation}) yields
\begin{equation*}
\frac{\big\arrowvert\E\big( {S}_n(\hat{x}_n,\hat{\delta}_n)-{S}_n(x,\delta_n)\big)\big\arrowvert}{\hat{\gamma}_n(\hat{x}_n,\hat{\delta}_n)}\ =\ O\Big(\delta_n^{-1}n^{-1/d}c_n\Big)\big(1+o_p(1)\big).
\end{equation*}
The arguments of the proof of Theorem \ref{thm: Effizienz2} apply again and lead to the expansion
\begin{align}
T_{\hat{j}_n\hat{k}_nn}\,-\, C_{\hat{j}_n\hat{k}_nn}\ &=\ O_p(1)\ +\ O\Big(\delta_n^{-1}n^{-1/d}c_n\Big)\big(1+o_p(1)\big)\nonumber\\  
&\ \ \ \ \ \ \ \ \ \ \ \ \ \,  +\ \frac{\E {S}_n(x,\delta_n)}{\gamma_n(x,\delta_n)}\Big(1\, +\, o_p(1)\Big)\ -\ \sqrt{2\log\Big({\gamma}_n(x,{\delta}_n)^{-2}\Big)},\label{ea}
\end{align}
while with the same reasoning as in the proof of Theorem \ref{thm: spatial adaptivity}
$$
\frac{\E {S}_n(x,\delta_n)}{\gamma_n(x,\delta_n)}\ \geq \ C\sqrt{n}\frac{c_n}{\sqrt{n\delta_n^d}}\delta_n^{d/2}\ =\ C c_n
$$
for some constant $C=C\big(d,(\lambda_n)_{\in\N}\big)>0$ and
$
\sqrt{2\log\big({\gamma}_n(x,{\delta}_n)^{-2}\big)}\ =\ O(1)\sqrt{\log(1/\delta_n)}
$. Thus, if $\sqrt{\log(1/\delta_n)}/c_n\rightarrow 0$ and $n\delta_n^d\rightarrow\infty$, (\ref{ea}) goes to infinity and the result follows.\hfill$\square$

\section{Appendix}
We start with a basic but useful property of the solution to (\ref{eq: rec 2}).
\begin{lemma}\label{lemm1}
If the solution to (\ref{eq: rec 2}) is not of bounded support, its lower isotonic and upper antitonic envelopes  are vanishing in $+\infty$. 
\end{lemma}

\vspace{-3.5mm}
\paragraph{\sc Proof} Suppose $\psi_{\beta}$ has only finitely many extremal points. From the last extremal point on the function $\psi_{\beta}$ is monotoneous and the integral over $\arrowvert .\arrowvert^{d-1}\psi_{\beta}(.)^2$ can only be finite if both envelopes are vanishing in $+\infty$. Now consider the case of infinitely many extremal points. 
Since the $L_2$-norm of the solution (\ref{eq: rec 2}) is finite and if there exists a sequence of local extrema of $\psi_{\beta}$ which stays uniformly bounded away from zero, their width must be bounded by a zero sequence. But now the result follows via contradiction of  (\ref{eq: function}), which, of course, is also applicable for local extrema.\hfill$\square$

\medskip
Let $\varepsilon>0$ be fixed. Define $t_{\varepsilon}$ to be a positive real number such that the following conditions are satisfied: $t_{\varepsilon}$ is a local extremal point, $\int_{B_0(t_{\varepsilon})}\gamma_{\beta}(x)^2d(x)\geq (1-\varepsilon/2)\Arrowvert\gamma_{\beta}\Arrowvert_2^2$, $\Arrowvert\psi_{\beta}\Arrowvert_{[t_{\varepsilon},\infty)}\leq \varepsilon/2$ (doable by Lemma \ref{lemm1}). Now extend the function $\psi_{\beta}I\{\cdot\leq t_{\varepsilon}\}$ to a compactly supported function $G_{\varepsilon}$ such that $G_{\varepsilon}(\Arrowvert .\Arrowvert_2) \in\HH_d(\beta,1;\R)$, $\int G_{\varepsilon}\big(\Arrowvert x\Arrowvert_2\big)dx=0$  and $\int_{\R^d\setminus B_{0}(t_{\varepsilon})}G_{\varepsilon}\big(\Arrowvert x\Arrowvert_2\big)^2dx$ smaller than $\varepsilon\Arrowvert \gamma_{\beta}\Arrowvert_2^2$; this is possible for $t_{\varepsilon}$ sufficiently large (because the uniform boundedness from $t_{\varepsilon}$ yields the boundedness of all partial derivatives by a multiple of $\varepsilon$ with the same argument as used in the proof of Theorem \ref{thm: spatial adaptivity}; so one may extend the function first to a compactly supported one in $\HH_d(\beta,L;\R)$ and then extend it close to zero such that its integral vanishes) - we omit an explicit construction at this point. With $\varepsilon$ sufficiently small, this construction leads to what was required in the proof of Theorem \ref{Thm: lower bound}.

\paragraph{\large\bf Acknowledgments} Lutz D\"umbgen's
contribution to the decoupling
subsequent to an extended discussion in Bern is gratefully acknowledged. Furthermore, I want to thank three unknown referees and an associate editor for their valuable comments and careful reading of the manuscript.

\begin{description}
\small
\item[]{\sc Behnen, K., Neuhaus, G. and Ruymgaart, F.} (1983).
       Two sample rank estimators of optimal nonparametric score-functions and
       corresponding adaptive rank statistics.\
       {\sl Ann.\ Statist.\ \textbf{11}}, 588--599.

\item[]{\sc Belomestny, D. and Spokoiny, V.} (2007).
       Spatial aggregation of local likelihood estimates with application to classification.\
       {\sl Ann.\ Statist.\ \textbf{35}}, 2287--2311.

\item[]{\sc Bennett, G.} (1962).\
       Probability inequalities for sums of independent random variables.\
       {\sl J.\ Amer.\ Statist.\ Assoc.\ \textbf{57}}, 33--45.

\item[]{\sc Butucea, C. and Tribouley, K.} (2006).\
       Nonparametric homogeneity tests.\ {\sl J.\ Statist.\ Plann. Inference\ \textbf{136}}, 597--639.

\item[]{\sc Donoho, D.} (1994a).
       Statistical estimation and optimal recovery.\
       {\sl Ann.\ Statist. \textbf{22}}, 238--270.

\item[]{\sc Donoho, D.} (1994b).
       Asymptotic minimax risk for sup-norm loss -- solution via optimal recovery.\
       {\sl Probab.\ Theory\ and\ Related \ Fields. \textbf{99}}, 145--170.

\item[]{\sc Ducharme, G.R. and Ledwina. T.} (2003).
       Efficient and adaptive nonparametric test for the two-sample problem.\
       {\sl Ann.\ Statist.\ \textbf{31}}, 2036--2058.

\item[]{\sc Dudley, R.M., Gin\'e, E. and Zinn, J.} (1991)
        Uniform and universal Glivenko-Cantelli classes.\ 
       {\sl J.\ Theoretical Probability.\ \textbf{4}}, 485--510.

\item[]{\sc D\"{u}mbgen, L.} (2002).
       Application of local rank tests to nonparametric regression.\
       {\sl J.\ Nonpar.\ Statist.\ \textbf{14}}, 511--537.


\item[]{\sc D\"{u}mbgen, L. and Spokoiny, V.G.} (2001).
	Multiscale testing of qualitative hypotheses.\
	{\sl Ann.\ Statist.\ \textbf{29}}, 124--152.

\item[]{\sc D\"{u}mbgen, L. and Walther, G.} (2008).
       Multiscale inference about a density.\
       {\sl Ann.\ Statist.\ \textbf{36}}, 1758--1785; {\sl accompagnying technical report available at} http://arxiv.org/abs/0706.3968

\item[]{\sc Eubank, R.L. and Hart, J.D.} (1992).
        Testing goodness-of-fit in regression via order selection criteria.\
        {\sl Ann.\ Statist.\ \textbf{20}}, 1412--1425.

\item[]{\sc Fan, J.} (1996).\
       Test of significance based on wavelet thresholding and Neyman's truncation.\
       {\sl J.\ Amer.\ Statist.\ Assoc.\ \textbf{91}}, 674--688.

\item[]{\sc H{\'a}jek, J. and $\check{\mathrm{S}}$idak, Z.} (1967).
       {\sl Theory of rank tests.} Academic press.

\item[]{Gijbels, I. and Heckmann, N.} (2004).
       Nonparametric testing for a monotone hazard function via normalized spacings.
       {\sl J.\ Nonpar.\ Statist.\ \textbf{16}}, 463--477. 

\item[]{\sc Hoeffding, W.} (1963).\
       Probability inequalities for sums of bounded random variables.\
       {\sl J.\ Amer.\ Statist.\ Assoc.\ \textbf{58}}, 13--30.

\item[]{Ingster, Y.} (1987).\
       Asymptotically minimax testing of nonparametric hypotheses.\
       {\sl Probability Theory and Mathematical Statistics, \textbf{I}}, 553--574.

\item[]{\sc Janic-Wr\'oblewska, A. and Ledwina, T.} (2000).
        Data driven rank test for two-sample problem.\ 
       {\sl Scand.\ J.\ Statist.\ \textbf{27}}, 281--297.

\item[]{\sc Klemel\"a, J. and Tsybakov, A.} (2001).
        Sharp adaptive estimation of linear functionals.\
        {\sl Ann.\ Statist.\ \textbf{29}}, 1567--1600.

\item[]{\sc Le Cam, L. and Yang, G.} (2000).
       {\sl Asymptotics in statistics: Some basic concepts.} Springer, New York.

\item[]{\sc Ledwina, T. and Kallenberg, W.C.M.} (1995).
        Consistency and Monte Carlo simulation of a data-driven version of smooth goodness-of-fit tests.\
        {\sl Ann.\ Statist.\ \textbf{23}}, 1594--1608.

\item[]{\sc Ledwina, T.} (1994).
        Data-driven version of Neyman's smooth test of fit.\
        {\sl J.\ Amer.\ Statist.\ Assoc.\ \textbf{89}}, 1000--1005.

\item[]{\sc Leonov, S.L.} (1997).
        On the solution of an optimal recovery problem and its applications in nonparametric statistics.\
        {\sl Math. Methods Statist.\ \textbf{4}}, 476--490.

\item[]{\sc Leonov, S.L.} (1999).
        Remarks on extremal problems in nonparametric curve estimation.\
        {\sl Statist. Probab. Lett.\ \textbf{43}}, 169--178.

\item[]{\sc Lepski, O. and Tsybakov, A.} (2000).  Asymptotically exact nonparametric hypothesis testing in sup-norm and at a fixed point. {\sl Probab.\ Theory\ Rel.\ Fields\ \textbf{117}}, 17--48.

\item[]{\sc Neuhaus, G.} (1982).
       $H_0$-contiguity in nonparametric testing problems and sample Pitman
       efficiency.\
       {\sl Ann.\ Statist.\ \textbf{10}}, 575--582. 

\item[]{\sc Neuhaus, G.} (1987).
        Local asymptotics for linear rank statistics with estimated score
        functions.\
        {\sl Ann.\ Statist.\ \textbf{15}}, 491--512.

\item[]{\sc Nussbaum, M.} (1996).
       Asymptotic equivalence of density estimation and Gaussian white noise.\
       {\sl Ann.\ Statist.\ \textbf{24}}, 2399-2430.

\item[]{\sc de la Pe$\tilde{\text{n}}$a, V.H.} (1994).
       A bound on the moment generating function of a sum of dependent variables with an application to simple sampling without replacement.\
       {\sl Ann.\ Inst.\ H.\ Poincar$\acute{\text{e}}$ \ Probab.\ Statist.\ \textbf{30}}, 197--211.

\item[]{\sc de la Pe$\tilde{\text{n}}$a, V.H.} (1999).
       A general class of exponential inequalities for martingales and ratios.\
       {\sl Ann.\ Prob.\ \textbf{27}}, 537--564.

\item[]{\sc Pollard, D.} (1984). {\sl Convergence of Stochastic Processes.}\ Springer.

\item[]{\sc Rohde, A.} (2008).
       Adaptive goodness-of-fit tests based on signed ranks.\
       {\sl Ann.\ Statist.\ \textbf{36}}, 1346--1374.

\item[]{\sc Rufibach, K. and Walther, G.} (2009).
       A block criterion for multiscale inference about a density, with applications to other multiscale problems.\
       {\sl J.\ Comp.\ Graph.\ Statist.}, to appear.

\item[]{\sc Serfling, R.J.} (1974).
       Probability inequalities for the sum of sampling without replacement.\
       {\sl Ann.\ Statist.\ \textbf{2}}, 39-48.

\item[]{\sc Shorack, G.R. and Wellner, J.A.} (1986).
       {\sl Empirical processes with applications to statistics.}\
       Wiley, New York.

\item[]{\sc Spokoiny, V.} (1996).
       Adaptive hypothesis testing using wavelets.\
       {\sl Ann.\ Statist.\ \textbf{24}}, 2477--2498.

\item[]{\sc van der Vaart, A.W. and Wellner, J.A.} (1996).
       {\sl Weak convergence and Empirical processes.}\
       Springer.

\item[]{\sc Walther, G.} (2010).
       Optimal and fast detection of spatial clusters with scan statistics.\
       {\sl Ann.\ Statist.\ \textbf{38}}, to appear.




\end{description}
\end{document}